\documentclass[12pt]{article}
\usepackage[fleqn]{amsmath}
\usepackage[cp1251]{inputenc}
\usepackage[russian,english]{babel}
\usepackage{amsmath,amsfonts,amssymb,graphicx,makeidx,amsthm,cite}
\usepackage{xcolor}
\usepackage{hyperref}
\usepackage[mathscr]{euscript}
 \let\mathscr\relax
\usepackage[scr]{rsfso}
\hypersetup{
	colorlinks = false,
	linkbordercolor = {violet},
}

\newtheorem{theorem}{Theorem}
\newtheorem{lemma}{Lemma}

\newtheorem{definition}{Definition}
\newtheorem{remark}{Remark}

\newtheorem{proposition}{Proposition}

\begin{document}
	\title{Fractional Wiener Chaos. Part 1.}
	
	\author{Elena Boguslavskaya\footnote{elena@boguslavsky.net , Department of Mathematics, Brunel University, Uxbridge UB8 3PH, United Kingdom}\, and Elina Shishkina\footnote{Department of Applied Mathematics,  Informatics  and Mechanics, Voronezh State University, Voronezh, Russia; Department of Applied Mathematics and Computer Modeling, Belgorod State National Research
			University (BelGU), Belgorod, Russia}}
	
	\maketitle
	
	
	\begin{abstract}
		{In this paper, we introduce a fractional analogue of the Wiener polynomial chaos expansion. It is important to highlight that the fractional order relates to the order of chaos decomposition elements, and not to the process itself, which remains the standard Wiener process.  The central instrument in our fractional analogue of the Wiener chaos expansion is the function denoted as $\mathcal{H}_\alpha(x,y)$, referred to herein as a power-normalised parabolic cylinder function.
			
			Through careful analysis of several fundamental deterministic and stochastic properties, we affirm that this function essentially
			serves as a fractional extension of the Hermite polynomial.
			In particular, the power-normalised parabolic cylinder function $\mathcal{H}_\alpha(W_t,t)$ demonstrates martingale properties and can be interpreted as a fractional It\^{o}  integral with 1 as the integrand,
			thereby drawing parallels with its non-fractional counterpart.
			
			To build a fractional analogue of polynomial Wiener chaos on the real line, we introduce a new function, which we call the extended Hermite function, by smoothly joining two power-normalized parabolic cylinder functions at zero. We form an orthogonal set of extended Hermite functions as a one-parameter family and use tensor products of the extended Hermite functions as building blocks in the fractional Wiener chaos expansion, in the same way that tensor products of Hermite polynomials are used as building blocks in the Wiener chaos polynomial expansion.
			
		} 
	\end{abstract} 
	\section{Introduction} 
	\label{sec:1}
	\setcounter{section}{1} \setcounter{equation}{0} 
	The area of fractional calculus has made its way into various pure and applied scientific fields,
	as evidenced by its integration into numerous disciplines, (see \cite{DietKiryLuchMachjTara2022}).
	Recently, there have been several attempts to incorporate fractional calculus into stochastic analysis (see e.g. \cite{Kolo2019,HachOkse2023}).
	In this paper, we create a fractional generalisation of  the Wiener chaos expansion (see \cite{DiNuOksePros2009,Engel1982,Schi2021}) by constructing an analogue of  Wiener chaos based on the extended Hermite functions.  We define these functions by smoothly joining at $0$ two specific parabolic cylinder functions with an exponential factor. 
	In the process, we demonstrate that this parabolic cylinder function with an exponential factor, which we call  ``a power-normalised parabolic cylinder function'', acts as an extension of a Hermite polynomial, particularly on $\mathbb{R^+}$. It retains the same martingale properties and other fundamental characteristics of a  Hermite polynomial. 
	However, power-normalised parabolic cylinder functions form an orthogonal basis only on  $\mathbb{R^+}$. To create an orthogonal basis on the entire $\mathbb{R}$, we need to smoothly join at $0$ two specific power-normalised cylinder functions. The resulting function, which we call the extended Hermite function, is used as a building block in the Wiener chaos expansion.
	Hence, it is accurate to state that power-normalised parabolic cylinder functions and extended Hermite functions can be viewed as fractional versions of Hermite polynomials.
	\\
	\\
	Many researchers have been exploring fractional Brownian motion and the associated stochastic analysis, as illustrated in books such as \cite{Mish2008} and \cite{BiagHuOkseZhan2008}. Here, we are pursuing another "fractional" aspect, focusing solely on ordinary Brownian motion while considering the order of chaos components to be fractional.
	\\
	\\
	The concept of homogeneous chaos, introduced by Wiener in 1938 (see \cite{Wiener1938}), involves studying square integrable nonlinear functionals of Brownian motion. This concept was reinterpreted by Cameron and Martin, who developed an orthogonal basis for these functionals using the so-called Fourier-Hermite functionals, formed from tensor products of Hermite polynomials. Later, in 1951, It\^{o} introduced the construction of multiple Wiener integrals to be used in the orthogonal decomposition of the space of square-integrable Brownian functionals. It is now well known, see e.g. \cite{Engel1982}, \cite{BiagHuOkseZhan2008}, that these three approaches describe the same concept. 
	\\
	In this paper, we propose a fractional extension of the polynomial chaos expansion. First, we introduce a fractional extension of the Hermite polynomial, which we call the power normalised cylinder function.  Using these power-normalised cylinder functions, we create an orthogonal system of functions on $\mathbb{R^+}$. Finally, by using two orthogonal sets of power normalised cylinder functions, we create the set of extended Hermite functions, which are orthogonal on $\mathbb{R}$.  We use the extended Hermite functions in the corresponding chaos , where the tensor product of Hermite polynomials in the  ordinary Wiener chaos expansion is replaced by the tensor product of extended Hermite functions.  This fractional analogue of Wiener chaos provides new expansions in orthogonal functions, which may be used in solving PDEs with such stochastic inputs when ordinary chaos expansion is not possible.
	\\
	\\
	Moreover,  the fractional analogue of the Hermite polynomial (power normalised cylinder function) can be considered  a fractional-order It\^{o} integral with integrand 1. The introduction of this fractional analogue of the Hermite polynomial is the small first step toward  developing the theory of fractional stochastic calculus. The question of whether the extended Hermite function can be considered as a a fractional-order It\^{o} integral, along with the other properties of this newly defined function,  will be discussed in the continuation of this paper (Fractional Wiener Chaos. Part 2 ), published separately.
  	\section{Preliminaries. Polynomial Wiener Chaos Expansion }
  	
	The Wiener Chaos expansion is a fundamental concept in stochastic analysis. Originating from the classical results of Wiener \cite{Wiener1938}, after the  work of It\^{o} \cite{Ito1951}, it now plays
	a central role in Malliavin calculus, see e.g. \cite{Nual2006}. The connection between multiple It\^{o} integrals and Hermite polynomials has become a classical result nowadays. In its polynomial form, the Wiener polynomial chaos is widely applied in computational engineering and physics \cite{GhanSpan1991}, \cite{Sudr2014}, \cite{Soiz2017}. The Wiener Chaos is also used in solving stochastic partial differential equations (SPDEs), when considered in the white noise framework, see e.g. \cite{HoldOkseUboeZhan2010}, \cite{LotoRozo2006b}.
	\\
	\\
	In this section, we present a brief introduction to Wiener polynomial chaos. For more detailed introduction  see  for example \cite{DiNuOksePros2009}, \cite{Jans1997}, \cite{GhanSpan1991}.
	\\
	\\
	The extended Hermite polynomials   (sometimes called Gaussian Hermite Polynomials), for $x \in \mathbb{R}$, $t>0$ and $n \in \mathbb{N} \cup {0}$  are  defined (see, for example, \cite{Schi2021}) as
	\begin{equation}
		\mathcal{H}_n(x,t) = (-1)^n t^n e^{\frac{x^2}{2t} }\frac{d^n}{dx^n}\left( e^{-\frac{ x^2}{2 t}}\right),\qquad \text{ for }
		n \in \mathbb{N}
		\label{ItoDefinitionHermitePolynomials}
	\end{equation}
	with $\mathcal{H}_0(x,t) =1.$ Sometimes, the extended Hermite polynomials are referred to as time-space Hermite polynomials, highlighting that $x$  can be considered as a spatial variable and
	$t$ as time. This perspective is particularly useful when dealing with
	$\mathcal{H}_n(W_t,t)$, where $(W_t)_{t \geq 0}$ is a Wiener process.
	\\
	\\
	For $t=1$ we get the classical Hermite polynomials $h_n(x) := \mathcal{H}_n(x,1)$, i.e.
	\begin{equation}
		h_n(x)=(-1)^ne^{\frac{x^2}{2}}\frac{d^n}{d x^n}e^{-\frac{x^2}{2}}.
	\end{equation}
	It is easy to see that
	\begin{equation}
		\label{hermit}
		\mathcal{H}_n(x,t)=t^{\frac{n}{2}}h_{n}\left(\frac{x}{\sqrt{t}}\right), \qquad n\in\mathbb{N}\cup\{0\}.
	\end{equation}
	The extended Hermite polynomials $\mathcal{H}_n(x,t)$ form an orthogonal basis in \\$L^2(\mathbb{R}, \nu(dx))$ with
	$\nu(dx) = \frac{1}{\sqrt{2 \pi t}}e^{-\frac{x^2}{2t}}dx$, see for example \cite{Scho2000}.
	\\
	\\
	Wiener \cite{Wiener1938}  originally introduced the homogeneous chaos expansion, utilising Hermite polynomials based on Gaussian random variables. Cameron and Martin's theorem \cite{CameMart1947} states that this expansion can approximate any functional within $L^2(C)$, achieving convergence in the sense of $L^2(C)$. This makes the Hermite chaos expansion a tool for representing second-order stochastic processes through orthogonal polynomials.
\vskip 0.2 cm
Let $g \in L^2 ([0,T])$ with
\begin{equation}
	||g||= \left(\int\limits_0^T g^2(t) dt\right)^{\frac{1}{2}},
\end{equation}
and let
\begin{equation} \xi= \int\limits_0^T g(t) dW_t,
\end{equation}
where $W=\{W_t\}_{t \geq 0}$is a Wiener process.
Note, that {$\frac{\Theta}{\|g\|}$} follows a standard normal distribution:
\begin{equation}
	\frac{\xi}{\|g\|}=\frac{1}{\|g\|}\int\limits_0^T g(t) dW_t \sim N(0,1).
\end{equation}
By formula (1.15) in \cite{DiNuOksePros2009}, the iterated It\^{o }integral can be expressed as follows for each natural $n$:
\begin{eqnarray}\label{RepeatedItoIntegralHermiteRelation}
	I_n (g)&=&n! \int\limits_0^T \int\limits_0^{t_n} \dots \int\limits_0^{t_2} g(t_1) g(t_2) \dots g(t_n) dW_{t_1} \dots dW_{t_n}= \|g\|^n h_n\left( \frac{\xi}{\|g\|} \right)  \nonumber\\
	&=& 
\mathcal{H}_n \left( \int\limits_0^T g(t) dW_t, \|g\|^2\right).
\end{eqnarray}
We take the latter expression as a  definition:
\begin{definition}
	\label{integral1}
	Let $g \in L^2([0,T])$ .
	We define the $n$-fold It\^{o} integral of function $g$ as:
	\begin{equation}
		I_n (g) :=  \mathcal{H}_n \left( \int\limits_0^T g(t) dW_t, \|g\|^2\right).
		\label{definition_I_n_g}
	\end{equation}
\end{definition}
We also make Proposition 1.8 from \cite{DiNuOksePros2009} into the following definition:
\begin{definition}
	\label{integral_tensorproduct}
	Let $g_1$, $g_2, \dots $ be orthogonal functions in $L^2([0,T])$.  We extend the definition \ref{integral1} to the tensor product of functions as follows:
	\begin{equation}
		I_n (g_{i_1}^{\otimes j_1}\!\! \otimes \dots \otimes  g_{i_r}^{\otimes j_r}) :=  \prod_{k=1}^{r}\mathcal{H}_{j_k} \!\!\left( \int\limits_0^T g_{i_k}(t) dW_t, \|g_{i_k}\|^2\right),
	\end{equation}
	where $\{j_k\}_{k=1}^r $is set of $r$ pairwise distinct natural numbers such that  $j_1 + \dots+ j_r =n$,  $\{i_k\}_{k=1}^r$ is also a set of $r$ pairwise distinct natural numbers, and  $\otimes$ denotes the tensor product.
\end{definition}
 Now consider the renormalisation of $\int_0^T g_k(t)dW_t$  such that the norms of functions $g_k$  are equal to 1, $\|g_k\|=1$ for all $k$.
Theorem 1.10 in \cite{DiNuOksePros2009} or Theorem 2.2.4 in \cite{HoldOkseUboeZhan2010}  
	 can be reformulated in our context as follows (see also \cite{GhanSpan1991}):
\begin{theorem}
	\label{WienerItoExpansion}
	Suppose $F$  is  a  square integrable random variable in $L^2(\mathbb{R}, \eta(x)dx)$ with $\eta(x)= \frac{1}{\sqrt{2 \pi }}e^{-\frac{x^2}{2}}$, and let  $\{ \xi_i\}_{i=1}^\infty$  be a set of independent normally distributed random variables with mean zero and variance one, $\xi_i \sim N(0,1)$. Then there exists a unique representation
	\begin{eqnarray}
		F &=& {\bf E} (F) + \sum_{n=1}^\infty\,\, \sum_{j_1\!+\dots\!+j_r\!=n} \,\,\sum_{\substack{i_1, i_2, \dots,i_r;\\ 
					 i_l \neq i_m\mbox{\small{for} } l\neq m}} c_{i_1 i_2 \dots i_r}^{j_1 j_2 \dots j_r} \prod_{k=1}^{r} h_{j_k}(\xi_{i_k}) \\&=& {\bf E} (F)+\sum_{n=1}^\infty \,\, \sum_{j_1\!+\dots\!+j_r\!=n} \,\,\sum_{\substack{i_1, i_2, \dots,i_r;\\ 
					 i_l \neq i_m\mbox{\small{for} } l\neq m}} c_{i_1 i_2 \dots i_r}^{j_1 j_2 \dots j_r} \prod_{k=1}^r \mathcal{H}_{j_k}(\xi_{i_k}, 1)
	\end{eqnarray}
	where $j_k $  and $i_k$ are natural numbers, $c_{i_1 i_2 \dots i_r}^{j_1 j_2 \dots j_r}$ are some constants. The convergence is in $L^2(\mathbb{R},\eta(x)dx)$ .
\end{theorem}

The represention  in Theorem \ref{WienerItoExpansion} is called the Wiener-It\^{o} expansion of random varialble $F$.

\section{Background and the Definition of the Power-Normalised Parabolic Cylinder Function}
 \label{sec:2}

        Let us consider the extended  Hermite polynomials depending on $x\in\mathbb{R}$ and some parameter $y>0$
\begin{equation}\label{Hn}
\mathcal{H}_n(x,y)=(-y)^ne^{\frac{x^2}{2y}}\frac{d^n}{dx^n}e^{-\frac{x^2}{2y}},\qquad n\in\mathbb{N}\cup \{0\}.
\end{equation}	
Formula \eqref{Hn} is called the {\it Rodrigues formula} for the (extended Hermite) polynomials.
\\
\\
Suppose the Gaussian density is given by
\begin{equation}
	\label{WDensity}
	w(x,y)={\frac {1}{\sqrt {2\pi y}}}e^{-\frac{x^{2}}{2y}},\qquad x \in \mathbb{R},\qquad y>0,
\end{equation}
	and $\mathbb{N}_{odd}$ and $\mathbb{N}_{even}$ denote the sets of all odd and even natural numbers, respectfully.
It is well known, see e.g. \cite{Foll2009}, that 
\begin{itemize}
	\item  $\{\mathcal{H}_n(x,y)\}_{n \in \mathbb{N}}$ forms an orthogonal basis on $\mathcal{L}^2(\mathbb{R}, w(x)dx)$,
	\item  $\{\mathcal{H}_n(x,y)\}_{n \in \mathbb{N}_{odd}}$ forms an orthogonal basis on $\mathcal{L}^2(\mathbb{R}^+, w(x)dx)$,
	\item  $\{\mathcal{H}_n(x,y)\}_{n \in \mathbb{N}_{even}}$ forms an orthogonal basis on $\mathcal{L}^2(\mathbb{R^+}, w(x)dx)$.
 \end{itemize}
In particular, 
\begin{equation}\label{ort}
	\int\limits_{-\infty}^\infty \mathcal{H}_n(\xi,y)\mathcal{H}_k(\xi,y)w(\xi,y)d\xi=n!y^{\frac{n+k}{2}}\delta_{nk},\qquad n,k\in\mathbb{N}\cup \{0\} 
\end{equation}
and
\begin{equation}\label{ort1}
\int\limits_{0}^\infty \mathcal{H}_n(\xi,y)\mathcal{H}_k(\xi,y)w(\xi,y)d\xi=n!y^{\frac{n+k}{2}}\delta_{nk},\qquad n,k\in\mathbb{N}_{even}\cup \{0\} \mbox{  or  } n,k\in\mathbb{N}_{odd}
\end{equation}
where $w=w(x,y)$ is
the Gaussian density (\ref{WDensity}),
$\delta_{nk}$ is the Kronecker delta. We aim to generalise the extended Hermite polynomials to fractional parameters while maintaining orthogonality.
\\
\\
In order to generalise extended Hermite polynomials to fractional parameters, consider the following representation: for  a non-negative integer $n$, we have (see \cite{AbraSteg1964}, formula 19.13.1)
\begin{equation}\label{H}
\mathcal{H}_n(x,y)=y^{\frac{n}{2}}e^{\frac{x^2}{4y}}D_n\left(\frac{x}{\sqrt{y}}\right),\qquad n\in\mathbb{N}\cup \{0\},
\end{equation} 	
where $D_n$
is the parabolic cylinder function.
\\
Now, for $\alpha\in\mathbb{C}$, the parabolic cylinder function $D_\alpha$
is given by
(see \cite{AbraSteg1964}, formula 19.12.1)
\begin{eqnarray}\label{parcyl}
&D_\alpha\left(z\right)=\sqrt{\pi}2^{\frac{\alpha}{2}}e^{-{\frac {z^{2}}{4}}}\times&\nonumber\\
\times &\left(\frac{1}{\Gamma \left({\frac {1-\alpha}{2}}\right)}
\,_{1}F_{1}\left(-{\frac{\alpha}{2}};{\frac {1}{2}};{\frac {z^{2}}{2}}\right){-}\frac{\sqrt {2}z}{\Gamma\left(-{\frac {\alpha}{2}}\right)}\,_{1}F_{1}\left(\frac {1-\alpha}{2};{\frac {3}{2}};{\frac {z^{2}}{2}}\right)\right).&
\end{eqnarray}	
Here $${}_{1}F_{1}(a;b;z)=\sum _{{n=0}}^{\infty }{\frac {(a)_nz^{n}}{(b)_n n!}},\qquad (a)_n=\frac{\Gamma(a+n)}{\Gamma(a)},
$$
is
a generalised hypergeometric function (a confluent hypergeometric  Kummer's function of the first kind).
Substituting in \eqref{parcyl} the expanded expressions of
hypergeometric series we get
\begin{equation}\label{SerD}
	D_\alpha\left(z\right){=}\frac{e^{-{\frac {z^{2}}{4}}}}{2^{\frac{\alpha}{2}+1}\Gamma(-\alpha)} \sum _{{n=0}}^{\infty }\frac {(-1)^n\Gamma\left({\frac{n-\alpha}{2}}\right)}{n!}(\sqrt{2}z)^{n}.
	\end{equation}
	Considering \eqref{H}, since  $\alpha\in \mathbb{C}$ in \eqref{parcyl}, we can
	extend $\mathcal{H}_n$  by the parameter $n$ to the set of all complex numbers.
	\begin{definition}
The function
\begin{equation} \label{Halpha01}
	\mathcal{H}_\alpha(x,y)=y^{\frac{\alpha}{2}}e^{\frac{x^2}{4y}}D_{\alpha}\left( \frac{x}{\sqrt{y}}\right),\qquad x\in\mathbb{R},\qquad y>0,\qquad \alpha\in \mathbb{C},
\end{equation}
is called  {\bf the power-normalised parabolic cylinder function}. Here, $D_{\alpha}$ denotes the parabolic cylinder function as described in \eqref{parcyl}.
\end{definition}
One can also consider $\mathcal{H}_\alpha(z,y)$ as a function of the complex variable $z$. Since
$D_{\alpha}(z)$ is an entire function of  $z$,  $\mathcal{H}_\alpha(z,y)$ is also an entire function of its argument $z$.
\section{Some Properties of the Power-Normalised Parabolic Cylinder Function}
\subsection{Partial Differential Equation and  Boundary Conditions for the Power Normalised Parabolic Cylinder Function }
\begin{proposition}\label{prop}
The power-normalised parabolic cylinder function $\mathcal{H}_\alpha(x,y)$, for $x,\alpha\in\mathbb{R}$, $y>0$,
has the following properties:
\begin{enumerate}
	\item For $x>0$ and as the parameter $y$ tends to zero,  $\mathcal{H}_\alpha(x,y)$ approached $x^\alpha$:
	\begin{equation}\label{Lim}
		\lim\limits_{y\rightarrow 0}	\mathcal{H}_\alpha(x,y)=x^\alpha.
	\end{equation}
	\item	Derivatives of $\mathcal{H}_\alpha(x,y)$ with respect to $x$ and  $y$, for $x\in\mathbb{R}$, $y>0$, are given by
	\begin{equation}\label{DxH}
		\frac{\partial}{\partial x}\mathcal{H}_\alpha(x,y)=\alpha\mathcal{H}_{\alpha-1}(x,y),\qquad
	\end{equation}	
	\begin{equation}\label{DtH}
		\frac{\partial}{\partial y}\mathcal{H}_\alpha(x,y)=-\frac{\alpha(\alpha-1)}{2}\mathcal{H}_{\alpha-2}(x,y).
	\end{equation}
	\item For $x>0$ and for any  $\alpha\in\mathbb{R}$ 	function $u(x,y)=\mathcal{H}_\alpha(x,y)$ is the solution $u=u(x,y)$ to the  problem:
	\begin{equation}
		\label{ZC01}
		\left\{ \begin{array}{l}
			u_y+\frac{1}{2}u_{xx}=0,\\
			u(x,0)= x^\alpha,\\
			u(0,y) = \sqrt{\pi}\frac{2^{\frac{\alpha}{2}} y^{\frac{\alpha}{2}}}{\Gamma(\frac{1 - \alpha}{2})}.
		\end{array} \right.
	\end{equation}
\end{enumerate}
\end{proposition}
\proof	1. According to the approximation $D_\alpha(z)\sim z^{\alpha}e^{-\frac{z^2}{4}}$ for large $z$  and moderate $\alpha$  as stated in  \cite{AbraSteg1964},
p. 689, formula 19.8.1 we get:
\begin{eqnarray}\label{Lim01}
\lim\limits_{y\rightarrow 0}	\mathcal{H}_\alpha(x,y)&=&\lim\limits_{y\rightarrow 0}y^{\frac{\alpha}{2}}e^{\frac{x^2}{4y}}D_\alpha\left(\frac{x}{\sqrt{y}}\right)\nonumber\\
&=&\lim\limits_{y\rightarrow 0}y^{\frac{\alpha}{2}}e^{\frac{x^2}{4y}}\left(\frac{x}{\sqrt{y}}\right)^{\alpha}e^{-\frac{x^2}{4y}}=x^\alpha.
\end{eqnarray}
2. For the derivative of $D_\alpha(z)$, we have the formula
\begin{equation} \label{DD}
\frac{d}{dz}D_\alpha(z)=\alpha D_{\alpha-1}(z)-\frac{z}{2}D_{\alpha }(z),
  \end{equation}
which can be obtained from 19.6.1 in \cite{AbraSteg1964}, p. 688.
Therefore, we can calculate
\begin{eqnarray*}
\frac{\partial}{\partial x}\mathcal{H}_\alpha(x,y)&=& y^{\frac{\alpha }{2}}\frac{\partial}{\partial x}\left[  e^{\frac{x^2}{4 y}}   D_{\alpha}\left(\frac{x}{\sqrt{y}}\right)\right]\\
&=&y^{\frac{\alpha }{2}}e^{\frac{x^2}{4 y}}\left[
\frac{x}{2y}D_\alpha\left(\frac{x}{\sqrt{y}}\right)+\alpha\frac{1}{\sqrt{y}}D_{\alpha-1}\left(\frac{x}{\sqrt{y}}\right)\right.
\\&-&
\left.	\left(\frac{x}{2\sqrt{y}}\right) \frac{1}{\sqrt{y}}D_\alpha\left(\frac{x}{\sqrt{y}}\right)\right] \\&=&
\alpha y^{\frac{\alpha-1 }{2}}e^{\frac{x^2}{4 y}}D_{\alpha-1}\left( \frac{x}{\sqrt{y}}\right)\\&=&\alpha\mathcal{H}_{\alpha-1}(x,y).
\end{eqnarray*}	
For the derivative with respect to $y$, utilising equation \eqref{DD} and formula 19.6.4 from \cite{AbraSteg1964}, p. 688, we get
\begin{eqnarray*}
\frac{\partial}{\partial y}\mathcal{H}_\alpha(x,y) &=& \frac{\partial}{\partial y}\left( y^{\frac{\alpha }{2}} e^{\frac{x^2}{4 y}}   D_{\alpha}\left(\frac{x}{\sqrt{y}}\right)\right)\\
& =&
\frac{\alpha }{2}y^{\frac{\alpha }{2}-1} e^{\frac{x^2}{4 y}} D_\alpha\left(\frac{x}{\sqrt{y}}\right)-\frac{x^2}{4}y^{\frac{\alpha }{2}-2} e^{\frac{x^2}{4 y}} D_\alpha\left(\frac{x}{\sqrt{y}}\right)\\& &+
\frac{x^2}{4}y^{\frac{\alpha }{2}-2}e^{\frac{x^2}{4 y}} D_\alpha\left(\frac{x}{\sqrt{y}}\right)-\frac{\alpha}{2}x y^{\frac{\alpha-3 }{2}} e^{\frac{x^2}{4 y}}D_{\alpha-1}\left(\frac{x}{\sqrt{y}}\right)\\&=&
\frac{\alpha }{2}y^{\frac{\alpha }{2}-1} e^{\frac{x^2}{4 y}}\left(  D_\alpha\left(\frac{x}{\sqrt{y}}\right)-\frac{x}{\sqrt{y}}D_{\alpha-1}\left(\frac{x}{\sqrt{y}}\right)\right) \\&=&
-\frac{\alpha(\alpha-1)}{2}y^{\frac{\alpha-2 }{2}} e^{\frac{x^2}{4 y}}D_{\alpha-2}\left( \frac{x}{\sqrt{y}}\right)\\ &=&- \frac{\alpha(\alpha-1)}{2}\mathcal{H}_{\alpha-2}(x,y).
\end{eqnarray*}
3. 	According to formulas \eqref{DxH} and \eqref{DtH} we get that $\mathcal{H}_\alpha(x,y)$ satisfies equation in \eqref{ZC01}.
According to formula \ref{Lim}  we get a condition in \ref{ZC01}. 

We   call $\mathcal{H}_\alpha(x,y)$  {\it the power-normalised parabolic cylinder function} due to the property \eqref{Lim}.
\subsection{The Norm of the Power Normalised Parabolic Cylinder Function}
\begin{figure}[ht] 
	\begin{center}
		\includegraphics[scale=0.4]{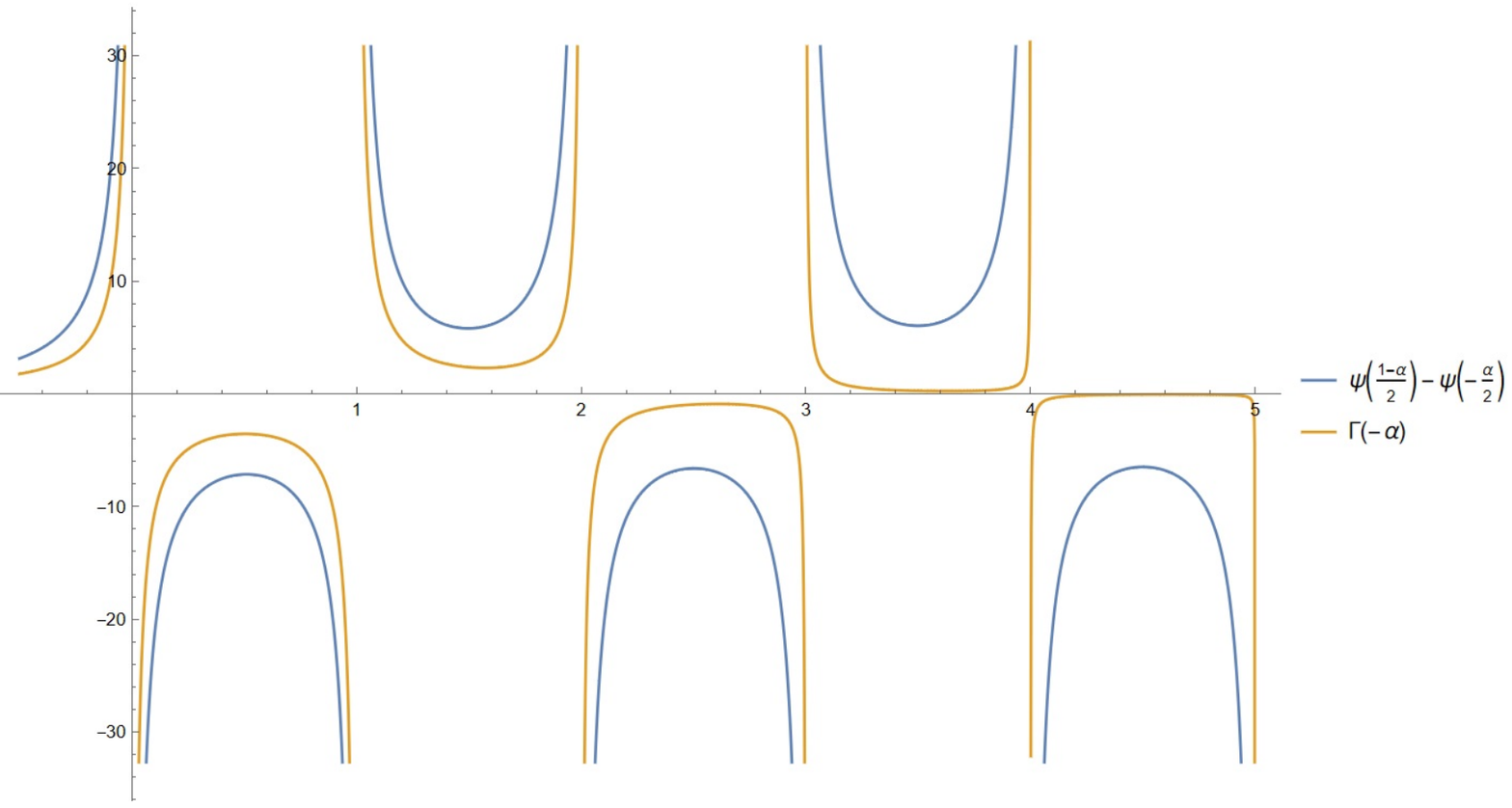}
		\caption{$\psi\left(\frac{1-\alpha}{2}\right)
			-\psi\left(-\frac{\alpha}{2}\right)$ and $\Gamma(-\alpha)$}\label{Plot01}
	\end{center}
\end{figure}

\begin{proposition}
	\label{PropNorm}
	For all real  $\alpha\neq 0,1,2,...$ and for fixed $y>0$
	\begin{equation}
		\label{N01}
		0<||\mathcal{H}_\alpha||_{w}^2=y^{\alpha} \frac{\psi\left(\frac{1-\alpha}{2}\right)
			-\psi\left(-\frac{\alpha}{2}\right)}{2\Gamma(-\alpha)},
	\end{equation}
	where $\psi$ is the digamma function, given by $\psi(z)=\frac{\Gamma'(z)}{\Gamma(z)}$.
	\end{proposition}
\proof We have (see \cite{GradRyzh2000}, 7.711.3)
\begin{eqnarray*}
||\mathcal{H}_\alpha||_{w}^2&=&\frac{1}{\sqrt{2\pi y}}\int\limits_0^\infty \mathcal{H}_\alpha^2(x,y)e^{-\frac{x^2}{2y}}dx\\
&=&\frac{1}{\sqrt{2\pi y}}y^{\alpha}\int\limits_0^\infty \left[e^{\frac{x^2}{4y}}D_{\alpha}\left( \frac{x}{\sqrt{y}}\right)\right]^2 e^{-\frac{x^2}{2y}}dx=y^{\alpha}\int\limits_0^\infty D_{\alpha}^2\left( \frac{x}{\sqrt{y}}\right)dx\\
&=&\frac{1}{\sqrt{2\pi y}}y^{\alpha}\int\limits_0^\infty D_{\alpha}^2\left( \frac{x}{\sqrt{y}}\right)dx\left|_{\left\{\frac{x}{\sqrt{y}}=\xi\right\}}\right.= \frac{1}{\sqrt{2\pi}}
y^{\alpha}\int\limits_0^\infty D_{\alpha}^2(\xi)d\xi\\
&=&y^{\alpha} \frac{\psi\left(\frac{1-\alpha}{2}\right)
	-\psi\left(-\frac{\alpha}{2}\right)}{2\Gamma(-\alpha)}.
\end{eqnarray*}
where $\psi$ is the digamma function, given by $\psi(z)=\frac{\Gamma'(z)}{\Gamma(z)}$.
\\
For $\alpha<0$ we get $\psi\left(\frac{1-\alpha}{2}\right)>\psi\left(-\frac{\alpha}{2}\right)$, since $\psi(z)$ is strictly increasing for $z>0$ and $\Gamma(-\alpha)>0$, therefore $||\mathcal{H}_\alpha||_{w}>0$.
\\
For $\alpha>0$, $\alpha\neq 0,1,2,...$ we can use the following equality (see \cite{Arfken2012}, p. 516)
$$
\psi(b+1)-\psi(a+1)=\sum\limits_{n=1}^{\infty}\frac{b-a}{(n+a)(n+b)}, \qquad a\neq b,\qquad a,b\neq -1,-2,-3,...
$$
Therefore,  we obtain the following formula
$$
\psi\left(\frac{1-\alpha}{2}\right)
-\psi\left(-\frac{\alpha}{2}\right)=\frac{1}{2}\sum\limits_{n=1}^{\infty}\frac{1}{\left(n-\frac{\alpha+2}{2}\right)\left(n-\frac{\alpha+1}{2}\right)},
$$
from which we see that   $\psi\left(\frac{1-\alpha}{2}\right)
-\psi\left(-\frac{\alpha}{2}\right)$ and $\Gamma(-\alpha)$
in \eqref{N01}  have the same simple poles at $\alpha=0,1,2,...$
and do not have any other poles.
Also, it is easy to see that $\psi\left(\frac{1-\alpha}{2}\right)
-\psi\left(-\frac{\alpha}{2}\right)\neq0$ and $\Gamma(-\alpha)\neq0$.
The behaviour of $\psi\left(\frac{1-\alpha}{2}\right)
-\psi\left(-\frac{\alpha}{2}\right)$ at $\alpha=k$, $k\in\mathbb{N}\cup\{0\}$ is given by (see Figure \ref{Plot01})
$$
\lim\limits_{\alpha\rightarrow 2k-0}\left(\psi\left(\frac{1-\alpha}{2}\right)
-\psi\left(-\frac{\alpha}{2}\right) \right)=
\lim\limits_{\alpha\rightarrow 2k-0}\Gamma(-\alpha)=+\infty,
$$
$$
\lim\limits_{\alpha\rightarrow 2k+0}\left(\psi\left(\frac{1-\alpha}{2}\right)
-\psi\left(-\frac{\alpha}{2}\right) \right)=
\lim\limits_{\alpha\rightarrow 2k+0}\Gamma(-\alpha)=-\infty,
$$
$$
\lim\limits_{\alpha\rightarrow 2k+1-0}\left(\psi\left(\frac{1-\alpha}{2}\right)
-\psi\left(-\frac{\alpha}{2}\right) \right)=
\lim\limits_{\alpha\rightarrow 2k+1-0}\Gamma(-\alpha)=-\infty,
$$
$$
\lim\limits_{\alpha\rightarrow 2k+1+0}\left(\psi\left(\frac{1-\alpha}{2}\right)
-\psi\left(-\frac{\alpha}{2}\right) \right)=
\lim\limits_{\alpha\rightarrow 2k+1+0}\Gamma(-\alpha)=+\infty.
$$
	Therefore, we have
	$$
	\frac{\psi\left(\frac{1-\alpha}{2}\right)
-\psi\left(-\frac{\alpha}{2}\right)}{2\sqrt{2}\Gamma(-\alpha)}>0.
$$
\subsection{ Rodrigues's Formula For the Power-Normalised  Parabolic Cylinder Function}
Here we present the fractional generalisation of Rodrigues's formula for extended Hermite polynomilas \ref{Hn}:
\begin{proposition} For $\mathcal{H}_\alpha(x,y)$ the following representation is valid
\begin{equation} \label{Halpha02}
\mathcal{H}_\alpha(x,y)=y^\alpha e^{\frac{x^2}{2y}}(\mathcal{D}_{-}^{\alpha})_x e^{-\frac{x^2}{2y}},\qquad \alpha>0,\qquad x \in \mathbb{R},\qquad y>0,
\end{equation}
where is $\mathcal{D}_{-}^{\alpha}$ is the Riemann-Liouville fractional derivative \eqref{RLD6}.
\end{proposition}
\begin{proof}
Let $m=[\alpha]+1$ when $\alpha$ in not integer and $m=\alpha$ when $\alpha$ is integer.
Calculating $(\mathcal{D}_{-}^{\alpha})_xe^{-\frac{x^2}{2y}}$ by  \eqref{RLD6} we obtain
\begin{eqnarray*}
(\mathcal{D}_{-}^{\alpha})_xe^{-\frac{x^2}{2y}}
&=&(-1)^m\frac{d^{m}}{dx^{m}}(I^{m-\alpha}_{-})_xe^{-\frac{x^2}{2y}}\\&=&
(-1)^my^{\frac{m-\alpha}{2}}\frac{d^{m}}{dx^{m}} e^{-\frac{x^2}{4y}}D_{m-\alpha}\left(\frac{x}{\sqrt{y}}\right).
\end{eqnarray*}
Here formula \eqref{LI1} when $\beta=m-\alpha$ was used.	
\\
By the formula (see \cite{BateErde1953}, p. 119. formula (16)) for $m\in\mathbb{N}$, $\alpha\in\mathbb{R}$
\begin{equation} \label{DU}
\frac{d^m}{dz^m}\left(e^{-\frac{z^2}{4}}D_a(z) \right)=(-1)^me^{-\frac{z^2}{4}}D_{a+m}(z)
\end{equation}
we obtain
$$
(\mathcal{D}_{-}^{\alpha})_xe^{-\frac{x^2}{2y}}=y^{-\frac{\alpha}{2}} e^{-\frac{x^2}{4y}}D_{\alpha}\left( \frac{x}{\sqrt{y}}\right).
$$
So, according to \ref{Halpha01} we get \ref{Halpha02}.
\end{proof}

The connection between parabolic cylinder function and fractional derivative was also noticed in \cite{PogaNada2021}.
\subsection{ On Zeros of  the Power-Normalised  Parabolic Cylinder Function}
It is known that  functions with infinitely many zeros are good candidates for constructing an orthogonal basis. Think, for example, of $\cos$ and $\sin$ functions. Let us show that  $\mathcal{H}_\alpha(z,y)$  fits this criterion. Specifically, we show that $\mathcal{H}_\alpha(z,y)$ has infinitely many zeros for each positive non-integer $\alpha$.  
Indeed, in  \cite{ElbMul1999}
the Hermite function
$$
{\bf H}_\alpha(x)=-\frac{\sin(\pi\alpha)\Gamma(\alpha+1)}{2\pi}\sum\limits_{n=0}^\infty\frac{\Gamma\left(\frac{n-\alpha}{2}\right)}{n!}(-2x)^n
$$
for $x\geq 0$ was considered.
From the  Euler's reflection formula,
we get
$$
\Gamma(-\alpha)=-\frac{\pi}{\sin(\alpha\pi)\Gamma(\alpha+1)}
$$
and, applying \eqref{SerD} and \eqref{Halpha01}, we obtain
$$
\mathcal{H}_\alpha(x,y)=\left(\frac{y}{2}\right)^{\frac{\alpha}{2}}{\bf H}_\alpha\left(\frac{x}{4\sqrt{2y}}\right).
$$
Therefore, from Theorem 3.1 in 
\cite{ElbMul1999}, we have the following proposition:
\begin{proposition}
	Let $x \geq 0$, $y>0$.
	For $n<\alpha\leq n+1$, $n=0,1,...,$ $\mathcal{H}_\alpha(x,y)$ has $n+1$ real zeros, and
	it has no real zeros when $\alpha\leq0$. Each zero is an increasing function of $\alpha$ on its interval
	of definition.
\end{proposition}

\section{Orthogonality of Power-Normalised Parabolic Cylinder Functions on The Half-Line}	
Here, we present the conditions under which
power-normalised parabolic cylinder functions
are orthogonal  with respect to the Gaussian weight $w(x) = \frac{1}{\sqrt{2\pi y}}e^{-\frac{x^2}{2y}}$ on $\mathbb{R}^+$.
\begin{theorem}\label{t00} Let $c$ be some real non-zero constant.
Consider the following equation: 
\begin{equation} \label{Eq01}
\frac{\Gamma \left(-\frac{\alpha}{2}\right)}{\Gamma \left(\frac{1-\alpha}{2}\right)}  = c ,
\qquad c \in\mathbb{R}\setminus{0}.
\end{equation}
Equation (\ref{Eq01} has infinitely many real positive non-integer roots. 
\\Suppose $\alpha_k>0$ and $\alpha_m>0$ with $k,m \in \mathbb{N}$ are real but not integer roots of (\ref{Eq01}).Then, the
functions $\{\mathcal{H}_{\alpha_k}(x,y)\}_{k\in\mathbb{N}}$ form an orthogonal set with respect to $x$ on the  interval $(0,\infty)$ for a fixed $y$, with
weight function $\frac{1}{\sqrt{2\pi y}}e^{-\frac{x^2}{2y}}$:
\begin{equation} \label{Eq02}
\langle \mathcal{H}_{\alpha_k},\mathcal{H}_{\alpha_m}\rangle_w=	\frac{1}{\sqrt{2\pi y}}\int\limits_0^\infty \mathcal{H}_{\alpha_k}(x,y)\mathcal{H}_{\alpha_m}(x,y)e^{-\frac{x^2}{2y}}dx=0,\qquad \alpha_k\neq\alpha_m,
\end{equation}
\begin{equation} \label{Eq03}
\frac{1}{\sqrt{2\pi y}}\int\limits_0^\infty \mathcal{H}_{\alpha_k}^2(x,y)e^{-\frac{x^2}{2y}}dx =y^{\alpha_k}\frac{\psi\left(\frac{1-\alpha_k }{2}\right)-\psi\left(-\frac{\alpha_k
	}{2}\right)}{2\Gamma (-\alpha_k )}.
	\end{equation}	
	
\end{theorem}
\proof  Let $\alpha_k$ and $\alpha_m$ be any real positive but not integer numbers.
For $\alpha_k\neq \alpha_m$ we consider the integral
\begin{equation*}
\int\limits_0^\infty \mathcal{H}_{\alpha_k}(x,y)\mathcal{H}_{\alpha_m}(x,y)e^{-\frac{x^2}{2y}}dx=y^{\frac{\alpha_k+\alpha_m}{2}}\int\limits_0^\infty D_{\alpha_k}\left( \frac{x}{\sqrt{y}}\right)D_{\alpha_m}\left( \frac{x}{\sqrt{y}}\right)dx.
\end{equation*}
Without loss of generality suppose $\alpha_k > \alpha_m$.
By formula 7.711.2 from \cite{GradRyzh2000} we have 
\begin{eqnarray}\label{OrthogonalZeros}
	&\int\limits_0^\infty D_{\alpha_k}\left( \frac{x}{\sqrt{y}}\right)D_{\alpha_m}\left( \frac{x}{\sqrt{y}}\right)dx \nonumber \\
	&= \frac{\pi 2^{\frac{1}{2}(\alpha_m+\alpha_k+1)}}{\alpha_k-\alpha_m}
	\left( \frac{1}{\Gamma(\frac{1-\alpha_k}{2})\Gamma(-\frac{\alpha_m}{2}) } - \frac{1}{\Gamma(\frac{1-\alpha_m}{2})\Gamma(-\frac{\alpha_k}{2}) }\right).
\end{eqnarray}
From the right side of  \eqref{OrthogonalZeros} we see  the if $\alpha_k$ and $\alpha_m$ are distinct positive non-integer roots of the equation (\ref{Eq01}), i.e. $\Gamma(-\frac{\alpha_m}{2}){=}c \Gamma(\frac{1 - \alpha_m}{2})$ and  $\Gamma(-\frac{\alpha_k}{2}){=} c\Gamma(\frac{1-\alpha_k}{2})$, then
\begin{eqnarray*}
& &	\frac{1}{\Gamma(\frac{1-\alpha_k}{2})\Gamma(-\frac{\alpha_m}{2}) } - \frac{1}{\Gamma(\frac{1-\alpha_m}{2})\Gamma(-\frac{\alpha_k}{2}) } \\& &=
	\frac{1}{c \Gamma(\frac{1-\alpha_k}{2})\Gamma(\frac{1-\alpha_m}{2}) } - \frac{1}{c \Gamma(\frac{1-\alpha_m}{2})\Gamma(\frac{1-\alpha_k}{2}) }=0.
\end{eqnarray*}
Moreover,  it is easy to see that if $\alpha_k$ and $\alpha_m$ are distinct positive natural numbers of the same parity, then the right-hand side in (\ref{OrthogonalZeros}) is also equal to zero, which proves the well-known fact that Hermite polynomials of the same parity form an orthogonal system in $\mathcal{L}^2(\mathbb{R}^+,w(x)dx)$.

\begin{figure}[ht] 
	\begin{center}
		\includegraphics[scale=0.7]{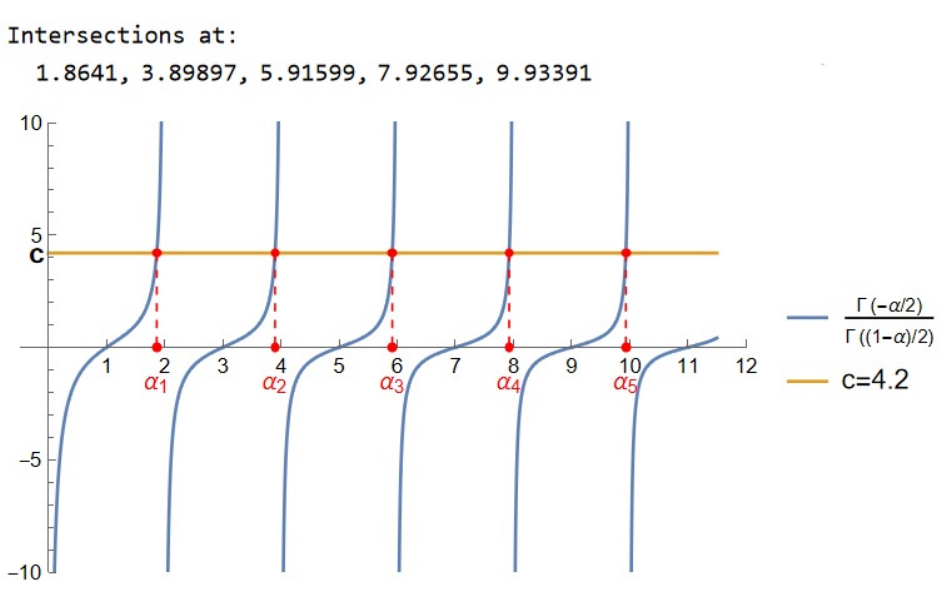}
		\caption{Roots of $A(\alpha)= c, $ with c=4.2. Here $\alpha_1 \approx 1.8641$, $\alpha_2 \approx 3.89897$, $\alpha_3 \approx 5.91599$, $\alpha_4 \approx 7.92655$, $\alpha_5 \approx 9.93391$. }\label{Ris02}
	\end{center}
\end{figure}

Now, let us show that (\ref{Eq01}) has infinitely many real positive not integer roots.  Denote
$$
A(\alpha)=\frac{\Gamma \left(-\frac{\alpha}{2}\right)}{\Gamma \left(\frac{1-\alpha}{2}\right)} 
$$
For $2n < \alpha < 2n+2$, $n\in \mathbb{N}\cup\{0\}$ the function $A=A(\alpha)$ is  continuous with $A(\alpha)<0$ for $2n< \alpha<2n+1$, $A(\alpha)>0$ for $2n+1< \alpha<2n+2$, and $A(2n+1)=0$. 
Moreover, we have for $n\in\mathbb{N} \cup 0$
$$
\lim\limits_{\alpha\rightarrow 2n+0} A(\alpha )=-\infty,
\qquad
\lim\limits_{\alpha\rightarrow 2n+2-0}A(\alpha)=+\infty.
$$
Therefore, equation $A(\alpha)=c$, where  $c \in \mathbb{R}\setminus{0}$ has infinitely many non-integer positive real roots, see Figure \ref{Ris02}.  (Similar function was investigated in \cite{Chadzitaskos2023})).
\\
The case $\alpha_k=\alpha_m$  has already been calculated in Proposition \ref{PropNorm} and is given by (\ref{N01}). 

\section{Sturm-Liouville Problem for the Hermite Equation and  Orthonormal Basis}
{In Theorem \ref{t00}, we discovered that the functions $\mathcal{H}_{\alpha_k}(x,y)$, where $k\in\mathbb{N}$ with $\alpha_k > 0$ and $\alpha_k \not\in {0,1,2,\dots}$ and $\frac{\Gamma\left(-\frac{\alpha_k}{2}\right)}{\Gamma\left(\frac{1-\alpha_k}{2}\right)}$=c, where $c$ is some fixed real non-zero number, form an orthogonal system with respect to $x$ in $L^2(\mathbb{R}^+,w(x)dx)$. In this section, we aim to address the question of the completeness of such a system. It is established (see \cite{Weidmann1980}) that identifying a complete orthogonal system of functions typically requires examining some boundary value problem for an ordinary differential equation, specifically, the Sturm-Liouville problem. The solution system of the specified boundary value problem then serves as a suitable basis. }\\\\
Let $y$ be a fixed parameter and $x\geq0$.
Now we would like to find the Sturm-Liouville problem for $\mathcal{H}_\alpha(x,y)$.
It is easy to see that
$$
x\frac{d}{dx}\mathcal{H}_\alpha(x,y)=\frac{x^2}{y}\mathcal{H}_\alpha(x,y)-\frac{x}{y}\mathcal{H}_{\alpha+1}(x,y),
$$	
$$
-y\frac{d^2}{dx^2}\mathcal{H}_\alpha(x,y)=-\frac{x^2}{y}\mathcal{H}_\alpha(x,y)+\alpha\mathcal{H}_\alpha(x,y)+\frac{x}{y}\mathcal{H}_{\alpha+1}(x,y),
$$	
therefore $\mathcal{H}_\alpha(x,y)$ satisfies to the  boundary value  problem
\begin{equation}\label{BVP01}
\left\{ \begin{array}{l}
L_\mathcal{H}\aleph(x)= \alpha \aleph(x),\qquad x>0;\\
\aleph(x)\underset{\substack{x\rightarrow+\infty}}{\sim}  x^\alpha,
\end{array} \right.
\end{equation}
where
\begin{equation*}
L_\mathcal{H}=-yD^2+xD,\qquad D=\frac{d}{dx}.
\end{equation*}
Here, we have taken into account \eqref{Lim01}.
\\
We also have
\begin{equation}
\label{xzero}	
\mathcal{H}_\alpha(0,y)=\frac{\sqrt{\pi } 2^{\frac{\alpha}{2}} y^{\frac{\alpha}{2}}}{\Gamma \left(\frac{1-\alpha }{2}\right)},\qquad
\frac{d}{dx}\mathcal{H}_\alpha(x,y)|_{x=0}=-\frac{\sqrt{\pi } 2^{\frac{\alpha +1}{2}} y^{\frac{\alpha-1}{2}}}{\Gamma \left(-\frac{\alpha }{2}\right)}.
\end{equation}	
So $\mathcal{H}_\alpha(x,y)$ is a nontrivial solution to the singular Sturm--Liouville problem \eqref{BVP01}.
\\
\\
{It is widely recognized that if an operator is self-adjoint within a specific domain, then several standard methods exist, see e.g. \cite{Weidmann1980}, to demonstrate that the eigenfunctions of this operator constitute an orthonormal basis in the corresponding Hilbert space. Consequently, to apply Sturm-Liouville theory to the operator  $L_\mathcal{H}$
, it is necessary to define a domain in which  $L_\mathcal{H}$
is self-adjoint for all functions within this domain. }
\\
Let $y>0$, $\xi>0$, $\xi\neq 0,1,2,...$, $a.c.(0,\infty)$ is the class of absolutely continuous on  $(0,\infty)$ functions,
$$\widetilde{D}=\{f,f'\in a.c.(0,\infty):f,L_\mathcal{H}f\in L^2(\mathbb{R}^+, w(x)dx)\},
$$
$a(\xi)=\Gamma\left(\frac{1-\xi}{2} \right)$, $b(\xi)=-\sqrt{\frac{y}{2}}\Gamma\left(-\frac{\xi}{2} \right)$.
Then in the definition domain
$$
D_{\xi}(L_\mathcal{H})=\left\{f\in\widetilde{D}:af(0)- bf'(0)=0\right\},\qquad a,b\in\mathbb{R},
$$
the operator $ L_\mathcal{H}$ is self-adjoint.
\\
Since
\begin{eqnarray}
& &\Gamma\left(\frac{1-\xi}{2} \right) \mathcal{H}_\alpha(0,y)+\sqrt{\frac{y}{2}}\Gamma\left(-\frac{\xi}{2} \right)\frac{d}{dx}\mathcal{H}_\alpha(0,y)|_{x=0}
\nonumber \\ \label{eq00}
&=&\frac{\sqrt{\pi } 2^{\frac{\alpha}{2}} y^{\frac{\alpha}{2}}}{\Gamma \left(\frac{1-\alpha }{2}\right)\Gamma \left(-\frac{\alpha }{2}\right)}
\left[\Gamma\left(\frac{1-\xi}{2} \right) \Gamma \left(-\frac{\alpha }{2}\right)-\Gamma\left(-\frac{\xi}{2} \right)\Gamma \left(\frac{1-\alpha }{2}\right)\right]\\&=&0 \nonumber,
\end{eqnarray}
we get that $\mathcal{H}_\alpha(x,y)\in D_{\xi}(L_\mathcal{H})$
by $x$.
\\
\\
{ In other words, $\mathcal{H}_\alpha(x,y)$ is an eigenfunction of the self-adjoint operator $L_\mathcal{H}$, with the corresponding parameter value $\alpha$ being its eigenvalue. Given that the problem is singular, the theory of eigenfunction expansion from the regular Sturm-Liouville problem cannot be applied directly. We shall consider the set of functions orthogonal with respect to $x$, $\{\mathcal{H}_{\alpha_k}(x,y)\}_{k\in\mathbb{N}}$, on the interval $(0,\infty)$ as described in Theorem \ref{t00}. Each $\mathcal{H}_{\alpha_k}(x,y)$ is an eigenfunction of $L_\mathcal{H}$, leading us to conclude that the eigenvalues $\alpha_k$, with $\alpha_k > 0$ and $\alpha_k \neq 1,2,\ldots$, are real, countable, ordered, and include a smallest eigenvalue. Hence, they can be listed as $\alpha_1 < \alpha_2 < \ldots < \alpha_k < \ldots$. However, there is no largest eigenvalue, and $\alpha_k$ approaches $+\infty$ as $k$ approaches $+\infty$. For each eigenvalue $\alpha_k$, there exists a corresponding eigenfunction $\mathcal{H}_{\alpha_k}(x,y)$. Thus, the operator $L_\mathcal{H}$ exhibits a purely discrete spectrum. According to \cite{Weidmann1980}, if a self-adjoint operator possesses a purely discrete spectrum, it must have a complete orthonormal sequence of eigenfunctions.
}
Therefore, eigenfunctions $\mathcal{H}_{\alpha_k}(x,y)$ corresponding to different eigenvalues form a basis in $L^2(\mathbb{R}^+,w(x)dx)$ with respect to $x$ under the conditions of Theorem \ref{t00}. Thus, we obtain the following theorem.
\begin{theorem}
For each $f(x) \in D_{\xi}(L_\mathcal{H})$, there exists a decomposition
\[
f(x) = \sum_{k=1}^{\infty} c_k \mathcal{H}_{\alpha_k}(x,1), \qquad c_k = \frac{\langle f(x), \mathcal{H}_{\alpha_k}(x,1) \rangle_w}{\|\mathcal{H}_{\alpha_k}\|_w},
\]
where $\alpha_k > 0$, $k \in \mathbb{N}$, and $\alpha_k \neq 0,1,2,\ldots$ are eigenvalues of $L_\mathcal{H}$ such that $\alpha_1 < \alpha_2 < \ldots < \alpha_k < \ldots$.
\end{theorem}
In \cite{Chadzitaskos2023}, an explicit form of a one-parameter family of
orthonormal bases for the space $L^2(\mathbb{R}^+,w(x)dx)$ was provided. Specifically, it demonstrated that the set of functions $D_\alpha$ forms an orthonormal basis in
$L^2(\mathbb{R}^+, w(x)dx)$ with eigenvalues $(\alpha+\frac{1}{2})$.

\section{Appell Integral Transform and the Martingale Property of the Power Normalised Parabolic Cylinder Function}
\subsection{Appell Integral Transform}
The main tool for studying the probabilistic properties of power-normalised parabolic cylinder functions of Wiener Process $\{\mathcal{H}_{\alpha}(W_t,t)\}_{t\geq 0}$ is the Appell integral transform, which was first introduced in (\cite{Bogu2014}) for the bilateral Laplace transform.  Here, we define the Appell integral transform on the Laplace transform to apply it to our purposes.
\\\\ We call $\mathcal{A}_{X_t} \{g\}(y)$ the Appell integral transform if
$$
\mathcal{A}_{X_t} \{g\}(y)=\mathcal {L}\, \frac{1}{\mathbf{E}e^{-yX_t}} \mathcal {L}^{-1}g,
$$
where $\mathcal {L}$ is the Laplace transform and $X_t$ is a stochastic process with $\mathbf{E}e^{-yX_t} < \infty$ for all $y \in \mathbb{R}$.  Given that we are dealing with generalised functions, it is important to discuss the appropriate class of functions for  $\mathcal{A}_{X_t}$.
Following \cite{Vladimirov1976}, we first introduce the action of the Fourier transform on the generalised function, and then we will narrow the scope of integration to the application of the Laplace transform.
\\
Let $f$ be a locally integrable function in $\mathbb{R}$ such that $f(u)=0$ for $u<0$ and $|f(u)|\leq Ae^{au}$ for $u\rightarrow+\infty$, then the Laplace transform of $f$ is
\begin{equation} \label{Lap}
\mathcal {L}\{f\}(s) =\int\limits_{0}^{\infty }e^{-su}f(u)\,du,\qquad s=\xi+i\eta\in\mathbb{C}.
\end{equation}	
Function $\mathcal {L}\{f\}(s)$ is analytic in the half-plane $\xi>a$.
\\
\\
The Fourier transform of an absolutely integrable on $\mathbb{R}$ function $f$
is
\begin{equation} \label{Four}
\mathcal {F}\{f\}(\xi)=\int\limits_{-\infty}^{\infty }e^{i\xi x}f(x)\,dx.
\end{equation}
Function $\mathcal {F}\{f\}(\xi)$ is continuous and bounded in $\mathbb{R}$.
\\
Suppose $\mathcal{D}=\mathcal{D}(\mathbb{R})$ is a  space of test function, $\mathcal{S}=\mathcal{S}(\mathbb{R})$ is a Schwartz space, and $\mathcal{D}'=\mathcal{D}'(\mathbb{R})$ and $\mathcal{S}'=\mathcal{S}'(\mathbb{R})$ t are their respective  dual spaces.
\\The following notation defines a canonical pairing:
$$\langle f,\varphi\rangle :=\int\limits_{-\infty}^\infty f(x)\varphi(x)\,dx,\qquad \varphi\in\mathcal{S}\qquad\text{or}\qquad\varphi\in\mathcal{D}.
$$
The Fourier transform  of any generalised   function $f\in \mathcal{S}'$ is given by
\begin{equation*}
\langle {\mathcal {F}}\{f\},\varphi \rangle=\langle f,{\mathcal {F}}\{\varphi\} \rangle,\qquad \varphi\in\mathcal{S}.
\end{equation*}	
The set of generalised functions $f\in\mathcal{D}'(\mathbb{R})$, $f=f(u)$, vanishing at $u<0$ is denoted by $\mathcal{D}'^+$. It is known  that if $f_1,f_2\in \mathcal{D}'^+$, then $f_1*f_2 \in \mathcal{D}'^+$ and $\mathcal{D}'^+$ is a convolution algebra. The unit element in the algebra $\mathcal{D}'^+$ is the $\delta$-function, since it satisfies $\delta*f=f$. Let $\mathcal{S}'^+=\mathcal{D}'^+\cap \mathcal{S}'$, $\mathcal{S}'^+$ be a convolution algebra.
Denote by $\mathcal{D}'^+_a$ the set of generalised functions $f\in \mathcal{D}'^+$, $f=f(u)$,   that have the property $f(u)e^{-\xi u}\in \mathcal{S}'^+$ for all $\xi>a$. $\mathcal{D}'^+_a$ is also a convolution algebra.
\\
\\
Let $f\in \mathcal{D}'^+_a$, $f=f(u)$, then $f(u)e^{-\xi u}\in \mathcal{S}'^+$ for all  $\xi >a$.  This implies that the generalised function  $f(u)e^{-\xi u}$ has a Fourier transform for each $\xi >a$, therefore the Laplace transform of the generalised function $f\in \mathcal{D}'^+_a$, $f=f(u)$, can be defined as
\begin{equation*}
\mathcal {L}\{f\}(s)=\mathcal{F}\{  f(u)e^{-\xi u}\}(-\eta)\in \mathcal{S}',\qquad s=\xi+i\eta\in\mathbb{C},\qquad \xi>a.
\end{equation*}	
Another way to define the action of Laplace transform $\mathcal{L}$ to   the distribution $f\in\mathcal{D}'^+_a$, $f=f(u)$,   for all $\xi>a$ is given by the formula
\begin{equation}\label{LapG}
\mathcal {L}\{f\}(s)=\langle e^{-au}f(u),\lambda(u)e^{-(s-a)u} \rangle,\qquad s=\xi+i\eta\in\mathbb{C},\qquad \xi>a,
\end{equation}	
where $\lambda(u)\in C^\infty$, $\lambda(u)=1$ in a neighborhood of the support of the function $f$ and $\lambda(u)=0$ for $u<-1$, therefore $\lambda(u)e^{-(s-a)u}\in\mathcal{S}$.
Formula \eqref{LapG} can be written in the form
\begin{equation*}
\mathcal {L}\{f\}(s)=\langle f(u),e^{-su}\rangle,\qquad s=\xi+i\eta\in\mathbb{C},\qquad \xi>a.
\end{equation*}	
Let $H_a$ be a set of functions $g=g(s)$ analytic in the half-plane $\xi>a$ such that for all $\varepsilon>0$ and $\xi_0>a$ there exist numbers $C_\varepsilon(\xi_0)\geq0$ and $m=m(\xi_0)\geq0$ that the inequality
\begin{equation}\label{NerLap}
|g(s)|\leq C_\varepsilon(\xi_0)e^{\varepsilon\xi}(1+|s|^m) ,\qquad \xi>\xi_0
\end{equation}	
is true. $H_a$ is an algebra with usual multiplication of analytical functions.
\\Let
$$
g(s)=\mathcal {L}\{f\}(s)=\langle f(u),e^{-su}\rangle,\qquad s=\xi+i\eta\in\mathbb{C},\qquad \xi>a.
$$
In order to $f\in \mathcal{D}'^+_a$,  it is necessary and sufficient that $g\in H_a$. Then for all $b\leq a$ and $\xi>\xi_0>a$
the inverse   Laplace transform of a function $g\in H_a$  is defined by
\begin{equation} \label{InvLap}
{\mathcal {L}}^{-1}\{g(s)\}(u)=f(u)={\frac {1}{2\pi i}}\left(\frac{d}{du}-b\right)^{m+2} \int\limits_{\xi -i\infty }^{\xi +i\infty }\frac{e^{us}g(s)}{(s-b)^{m+2}}\,ds,
\end{equation}	
where $m$ is from \eqref{NerLap}, gives $f\in \mathcal{D}'^+_a$. Thus, a one-to-one correspondence is established between the algebras $\mathcal{D}'^+_a$ and $H_a$.
\\
\\
Suppose $f(u)={\mathcal {L}}^{-1}\{g\}(u)$ is a generalised function and  $\{X_t\}_{t \geq 0}$ is a stochastic process such that $\frac{f(u)}{\mathbf{E}e^{-uX_t}}\in \mathcal{D}'^+_a$.
Consider the pairing
\begin{equation}\label{Appelltransform01}
\left\langle  \frac{f(u)}{\mathbf{E}e^{-uX_t}}, e^{-uy} \right\rangle=\left\langle   f(u), \frac{e^{-uy}}{\mathbf{E}e^{-uX_t}} \right\rangle\in H_a,
\end{equation}
with respect to $X=\{X_t\}_{t \geq 0}$.
For $g(s)=\mathcal {L}\{g\}(s)\in H_a$ we can rewrite \eqref{Appelltransform01} in the form of an integral transform
\begin{equation}\label{Appelltransform02}
\mathcal{A}_{X_t} \{g\}(y)=\left\langle  {\mathcal {L}}^{-1}g(u), \frac{e^{-uy}}{\mathbf{E}e^{-uX_t}} \right\rangle.
\end{equation}
We can formally  express $\mathcal{A}_{X_t}$ as
\begin{equation}\label{Appelltransform04}
\mathcal{A}_{X_t} \{g\}(y)=\mathcal {L}\, \frac{1}{\mathbf{E}e^{-uX_t}} \mathcal {L}^{-1}g
\end{equation}
which we understand in the sense of \eqref{Appelltransform02}.
\\
\\
For convenience, whether $\mathcal {L}^{-1}g(u)$ is a regular functional or a singular functional, we will use the same notation
\begin{equation}\label{Appelltransform03}
\mathcal{A}_{X_t} \{g\}(y)=\int\limits_{0}^\infty {\mathcal {L}}^{-1}\{g(s)\}(u)\frac{e^{-uy}}{\mathbf{E}e^{-uX_t}}du
\end{equation}
and understand (\ref{Appelltransform03}) in the sense \eqref{Appelltransform01}.
\\
We will refer to the transform given in equation \eqref{Appelltransform03} as the 'Appell Integral Transform'.
\begin{theorem}\label{t01}  If
\begin{enumerate}
\item $\mathcal {L}^{-1}g\in \mathcal{S}'$ is a generalised function
such that $\frac{\mathcal {L}^{-1}g}{\mathbf{E}e^{-uX_t}}\in \mathcal{D}'^+_a$,
\item $\varphi(y,u)=\frac{e^{-uy}}{\mathbf{E}e^{-uX_t}}\rho_{X_t}(y)\in\mathcal{S}(\mathbb{R}\times\mathbb{R}_+)$, where $\rho_{X_t}(y)$ is the density of $X_t$,
\item If for each fixed $y$, $y>0$ $\frac{e^{-yX_t}}{\mathbf{E}e^{-yX_t}}$ is a martingale with respect to $\sigma$-algebra $\mathcal{F} = (\mathcal{F}_t)_{t \geq 0}$,
\end{enumerate}
then $ \mathcal{A}_{X_t} \{g\}(X_t)$ is also a martingale with respect to $\sigma$-algebra $\mathcal{F} $.
\end{theorem}
\proof
Firstly,  we see that  $\mathbf{E} |\mathcal{A}_{X_t}\{g\}(X_t)| < \infty$.
\\
Secondly, in order to show that $\mathcal{A}_{X_t}\{g\}(X_t)$ is a martingale we should show that
\begin{equation*}
\mathbf{E} [\mathcal{A}_{X_t} \{g\}(X_t)|\mathcal{F}_s]= \mathcal{A}_{X_s} \{g\}(X_s).
\end{equation*}	
Let $\rho_{X_t|\mathcal{F}_s}(y)$ denote conditional density of $X_t$ with respect to sigma algebra $\mathcal{F}_s$.  Then for the conditional mathematical expectation with respect $\mathcal{F}_s$ we can write
\begin{eqnarray*}
\mathbf{E} [\mathcal{A}_{X_t} \{g\}(X_t)|\mathcal{F}_s]
&=&\mathbf{E}\left[\left. \int\limits_0^\infty   \mathcal {L}^{-1}g(u)\frac{e^{-uX_t}}{\mathbf{E}e^{-uX_t}}du \right| \mathcal{F}_s \right]\\
&=&\int\limits_{-\infty}^\infty  \left( \int\limits_{0}^\infty
\mathcal {L}^{-1}g(u)\frac{e^{-uy}}{\mathbf{E}e^{-uX_t}}du\right) \rho_{X_t|\mathcal{F}_s}(y)\,dy\\&=& \int\limits_{-\infty}^\infty  \left\langle \mathcal {L}^{-1}g(u),\frac{e^{-uy}}{\mathbf{E}e^{-uX_t}} \rho_{X_t|\mathcal{F}_s}(y) \right\rangle \,dy.
\end{eqnarray*}
Since $\mathcal {L}^{-1}g\in \mathcal{S}'$ and $\frac{e^{-uy}}{\mathbf{E}e^{-uX_t}}\rho_{X_t}(y)\in\mathcal{S}(\mathbb{R}\times\mathbb{R}_+)$ we can write
\begin{eqnarray*}
\mathbf{E} [\mathcal{lA}_{X_t} \{g\}(X_t)|\mathcal{F}_s]
&=&  \left\langle \mathcal {L}^{-1}g(u),\int\limits_{-\infty}^\infty  \frac{e^{-uy}}{\mathbf{E}e^{-uX_t}} \rho_{X_t|\mathcal{F}_s}(y) \,dy\right\rangle \\
&=&  \left\langle \mathcal {L}^{-1}g(u),\mathbf{E}\left.\left[\frac{e^{-uX_t}}{\mathbf{E}e^{-uX_t}} \right| \mathcal{F}_s\right] \right\rangle\\
&=&  \left\langle \mathcal {L}^{-1}g(u), \frac{e^{-uX_s}}{\mathbf{E}e^{-uX_s}}  \right\rangle =\mathcal{A}_{X_s} \{g\}(X_s),
\end{eqnarray*}
as $\mathbf{E}\left.\left[\frac{e^{-uX_t}}{\mathbf{E}e^{-uX_t}} \right| \mathcal{F}_s\right] = \frac{e^{-uX_s}}{\mathbf{E}e^{-uX_s}}  $ follows from the assumption  that $\frac{e^{-u X_t}}{\mathbf{E}e^{-u X_t}}  $ is a martingale with respect to $\mathcal{F}$.

\subsection{Power-Normalised Parabolic Cylinder Function as an Appell Integral transform}
Let $W=\{W_t\}_{t \geq 0}$ be a Wiener process. For $t>0$ the density function of this process is given by $w(x,t)$ (see \eqref{WDensity}).
Formally, we agree that $w(x,0)=\delta(x)$. Therefore, function $w(x,0)$ is a fundamental solution to the heat equation
\begin{equation*}
\frac{\partial w}{\partial t} =\frac{1}{2}	\frac{\partial^2 w}{\partial x^2}.
\end{equation*}
Therefore,    for $t\geq0$ we can write the Markov operator $P_t$ as
\begin{equation}\label{ex01}
P_tg(x)=\int\limits_{-\infty}^\infty g(s)w(x-s,t)ds={\frac {1}{\sqrt {2\pi t}}}\int\limits_{-\infty}^\infty g(y)e^{-\frac{(x-y)^2}{2t}}dy.
\end{equation}
Now, since
\begin{eqnarray}\label{mo}
\mathbf{E}e^{-u W_{t}}&=\int\limits_{-\infty }^{\infty }e^{-ux}w(x,t)dx=\frac {1}{\sqrt {2\pi t }}\mathcal{L}\{e^{-\frac {x^2}{2t}}\}(u)\nonumber\\
&=\frac {1}{\sqrt {2\pi t }}\int\limits_{-\infty }^{\infty }e^{-u x-\frac {x^2}{2t}}dx=e^{\frac{u^2}{2}t},
\end{eqnarray}
then the Appell integral transform  \eqref{Appelltransform03} for Wiener process is
\begin{equation}\label{Appelltransform0W0}
\mathcal{A}_{W_t} \{g\}(y)=\left\langle  {\mathcal {L}}^{-1}g(u), e^{-uy-\frac{u^2}{2}t}\right\rangle,\qquad g\in \mathcal{S}'
\end{equation}
or
\begin{equation}\label{Appelltransform0W}
\mathcal{A}_{W_t} \{g\}(y)=\int\limits_{0}^\infty  \mathcal{L}^{-1}g(u)e^{-uy-\frac{u^2}{2}t}du,\qquad g\in \mathcal{S}'
\end{equation}
or
\begin{equation}\label{Appelltransform0W01}
\mathcal{A}_{W_t} \{g\}(y)=\mathcal {L}\, e^{-\frac{u^2}{2}t} \mathcal{L}^{-1}g,\qquad g\in \mathcal{S}'.
\end{equation}
\begin{theorem}\label{t02}
For $\alpha\in\mathbb{R}$ function $\mathcal{H}_\alpha(x,t)$ can be represented as the Appell integral transform for Wiener process in the form
\begin{equation} \label{MT}
\mathcal{H}_{\alpha}(x,t)=\mathcal{A}_{W_t} \{y^{\alpha}\}(x).
\end{equation}
\end{theorem}
\proof
Let  first consider the case $\alpha<0$. Denote $\beta=-\alpha>0$.
According to  formula 1.4.58 on p.23 in \cite{KilbSrivTruj2006},  we have
$$
{\mathcal {L}}^{-1}\{y^{-\beta}\}(u)=\frac{u^{\beta-1}}{\Gamma(\beta)},\qquad {u>0}.
$$
Employing \eqref{Appelltransform0W},	we obtain
\begin{eqnarray*}
\mathcal{A}_{W_t} \{y^{-\beta}\}(x)&=&\int\limits_{0}^\infty \mathcal{L}^{-1}\{y^{-\beta}\}(u)e^{-ux-\frac {u^2}{2}t}du
\\&=&
\frac{1}{\Gamma(\beta)}\int\limits_{0}^\infty u^{\beta-1}e^{-ux-\frac {u^2}{2}t}du \\&=&
\frac{1}{\Gamma(\beta)}e^{\frac{x^2}{2t}}\cdot e^{-\frac{x^2}{2t}}\int\limits_{0}^\infty u^{\beta-1}e^{-ux-\frac {u^2}{2}t}du
\\&=&
\frac{1}{\Gamma(\beta)}e^{\frac{x^2}{2t}}\int\limits_{0}^\infty u^{\beta-1}e^{-\frac{1}{2t}(ut+x)^2}du\left|_{\{ut+x=z\}}\right.\\&=&
\frac{1}{t^\beta\Gamma(\beta)}e^{\frac{x^2}{2t}}\int\limits_{x}^\infty (z-x)^{\beta-1}e^{-\frac{z^2}{2t}}dz=\frac{1}{t^\beta}e^{\frac{x^2}{2t}} ((I^\beta_{-})_ze^{-\frac{z^2}{2t}})(x).
\end{eqnarray*}
Applying formula \eqref{LI1} and \eqref{Halpha01}  we can write
$$
\mathcal{A}_{W_t} \{y^{-\beta}\}(x)=  t^{-\frac{\beta}{2}}e^{\frac{x^2}{4t}}D_{-\beta}\left(\frac{x}{\sqrt{t}}\right)=\mathcal{H}_{-\beta}(x,t).
$$
This results in \eqref{MT} for $\alpha<0$.
\\
\\
Now, let us prove \eqref{MT} for $\alpha\geq0$. We should consider two scenarios: when $\alpha$ is an integer and when $\alpha$ is not an integer.
\\
Let $\alpha=n\in \mathbb{N}\cup\{0\}$, then
the inverse   Laplace transform of $y^n$ is the $n$-th derivative of the delta
function (see \cite{KilbSrivTruj2006}, p. 24, formula 1.4.62)
$$
{\mathcal {L}}^{-1}\{y^n\}(u)=\delta^{(n)}(u).
$$
Since 	by \eqref{hermit} $\frac{d^n}{du^n} \left( e^{-ux-{\frac {u^2}{2}}t}\right) _{u=0}=(-1)^n\mathcal{H}_n(x,t)$, then by
\eqref{Appelltransform0W0} we get
\begin{eqnarray*}
\mathcal{A}_{W_t}\{y^{n}\}(x) &=&\left\langle  \delta^{(n)}(u), e^{-uy-\frac{u^2}{2}t}\right\rangle\\
&=&(-1)^n \left. \left( \frac{d^n}{du^n} e^{-ux-{\frac {u^2}{2}}t}\right)\right|_{u=0}\\&=&\mathcal{H}_n(x,t).
\end{eqnarray*}
Finally, let $\alpha>0$ and not integer.
For such $\alpha$ and for the Caputo fractional derivative \eqref{RLDC} we have (see \cite{FengYeZhang2020})
$$
{\mathcal {L}}^{-1}\{y^\alpha\}(u)=	({}^C\mathcal{D}_{0+}^{\alpha}\delta)(x).
$$
Here we used Caputo derivative for convenience, because next we apply the formula for integration by parts  in fractional integrals and we should obtain Riemann-Liouville integral which is connected with $\mathcal{H}_\alpha$ (see \eqref{Halpha02}).
\\
\\
Let $m=[\alpha]+1$. Then by   \eqref{Appelltransform0W0} and by formula of integrating by parts in fractional integrals (see \cite{SamkKilbMari1993}) we obtain
\begin{eqnarray*}
\mathcal{A}_{W_t}\{y^{\alpha}\}(x)&=&\left\langle (\,^C\mathcal{D}_{0+}^{\alpha}\delta)(u) , e^{-uy-\frac{u^2}{2}t}\right\rangle\\&=&\left\langle  \delta^{(m)}(u), ((I^{m-\alpha}_{-})_u e^{-ux-\frac{u^2}{2}t} \right\rangle.
\end{eqnarray*}
By \eqref{LI2} we get
$$
((I^{m-\alpha}_{-})_u e^{-xu-\frac{u^2}{2}t} )(\tau)=t^{-\frac{m-\alpha}{2}}	e^{\frac{x^2}{2t}-\frac{(x+t\tau)^2}{4t}}D_{\alpha-m}\left(\frac{x+t\tau}{\sqrt{t}}\right).
$$
According to expansion of $D_\beta(z)\sim z^{\beta}e^{-\frac{z^2}{4}}$ for $z$ large and $\beta$ moderate from \cite{AbraSteg1964},
p. 689, formula 19.8.1, function $(I^{m-\alpha}_{-})_u e^{-ux-\frac{u^2}{2}t} $ can be considered as a Schwartz  test function by $u$ and
$$
\mathcal{A}_{W_t}\{y^{\alpha}\}(x)=\left\langle  \delta^{(m)}(u), t^{\frac{\alpha-m}{2}}	e^{\frac{x^2}{2t}-\frac{(x+tu)^2}{4t}}D_{\alpha-m}\left(\frac{x+tu}{\sqrt{t}}\right)\right\rangle.
$$
Using the representation of parabolic cylinder function $U(\beta,z)$ as
an integral along the real line for $\gamma>-\frac{1}{2}$ (see \cite{AbraSteg1964}, p. 687, formula 19.5.3)
\begin{equation} \label{U}
D_{-\gamma-\frac{1}{2}}(z)=
\frac{e^{\frac{z^2}{4}}}{\Gamma\left(\gamma+\frac{1}{2}\right) }\int\limits_0^{\infty}\tau^{\gamma-\frac{1}{2}}e^{-\frac{(\tau+z)^2}{2}}d\tau
\end{equation}
and \eqref{Halpha01}, we get
\begin{eqnarray*}
(-1)^m\mathcal{A}_{W_t}\{y^{\alpha}\}(x)&=&\left.   t^{\frac{\alpha-m}{2}} e^{\frac{x^2}{2t}} \frac{d^{m}}{du^{m}} e^{-\frac{(x+tu)^2}{4t}}D_{\alpha-m}\left(\frac{x+tu}{\sqrt{t}}\right)\right|_{u=0}\\
&  =&\left. (-1)^m e^{\frac{x^2}{2t}} t^{\frac{\alpha}{2}} e^{-\frac{(x+tu)^2}{4t}}D_\alpha\left(\frac{x+tu}{\sqrt{t}}\right)\right|_{u=0}
\\&=& (-1)^m e^{\frac{x^2}{4t}} t^{\frac{\alpha}{2}} D_{\alpha}\left( \frac{x}{\sqrt{t}}\right)=\\&=&(-1)^m	\mathcal{H}_\alpha(x,t).
\end{eqnarray*}
and
$\mathcal{A}_{W_t}\{y^{\alpha}\}(x)=\mathcal{H}_\alpha(x,t)$.
That yields \eqref{MT} for non-integer $\alpha>0$.

\begin{remark}
It is worth noting that, similar to Theorem \ref{MT}, it can be demonstrated that
\begin{equation}
\mathcal{A}_{\Theta_t}\{y^\alpha\}(\Theta_t) = \mathcal{H}_\alpha (\Theta_t, {\bf E} \Theta_t^2)  =\mathcal{H}_\alpha (\Theta_t, \|g\|_{L^2}),
\end{equation}
where $\Theta_t =\int_0^t g(s) dW_s$ with $g \in L^2$ and $(W_t)_{t \geq 0}$ is a Wiener process.
\end{remark}
\section{Power-normalised Parabolic Cylinder Functions as a Stochastic Process and Its Properties}
Let $\{W_t\}_{t\geq 0}$ be a Wiener process,  $\alpha>0$, then  $\mathcal{H}_\alpha(W_t,t)$ is a stochastic process.
We shall provide some probabilistic properties of the process $\mathcal{H}_\alpha(W_t,t)$.
\begin{lemma} Let $\{W_t\}_{t\geq 0}$ be a Wiener process. Then the following formulas are valid
\begin{eqnarray*}
d\mathcal{H}_\alpha(W_t,t)&=&	 \alpha\mathcal{H}_{\alpha-1}(W_t,t)dW_t,\\
\int\limits_0^t \mathcal{H}_{\alpha-1}(W_s,s))dW_s&=&\frac{1}{\alpha}\mathcal{H}_\alpha(W_t,t).
\end{eqnarray*}
\end{lemma}
\proof
We know from Definition \ref{Halpha01} that $\mathcal{H}_\alpha(x,y)$ has all derivatives by $x$ and by $t>0$. For $t=0$ we take a limit with $t\to+0$ of $\mathcal{H}_\alpha$ or its derivatives by $t$. Therefore, we can apply  It\^{o} formula
\begin{eqnarray*}
d\mathcal{H}_\alpha(W_t,t)&=&\\&=&\frac{\partial \mathcal{H}_\alpha}{\partial x}(W_t,t)dW_t+
\frac{1}{2}\frac{\partial^2 \mathcal{H}_\alpha}{\partial x^2}(W_t,t)(dW_t)^2+\frac{\partial \mathcal{H}_\alpha}{\partial t}(W_t,t)dt\\
&=&\frac{\partial \mathcal{H}_\alpha}{\partial x}(W_t,t)dW_t+
\left(\frac{\partial \mathcal{H}_\alpha}{\partial t}(W_t,t)+\frac{1}{2}\frac{\partial^2 \mathcal{H}_\alpha}{\partial x^2}(W_t,t)\right) dt\\&=&\frac{\partial \mathcal{H}_\alpha}{\partial x}(W_t,t)dW_t\\
&=& \alpha\mathcal{H}_{\alpha-1}(W_t,t))dW_t
\end{eqnarray*}
we have $$d\mathcal{H}_\alpha(W_t,t)=	 \alpha\mathcal{H}_{\alpha-1}(W_t,t))dW_t$$ and
$$
\int\limits_0^t\mathcal{H}_{\alpha-1}(W_s,s))dW_s=\frac{1}{\alpha}\mathcal{H}_{\alpha}(W_t,t).
$$

\begin{theorem}
If $\{W_t\}_{t\geq 0}$ is a Wiener process, then the process $\mathcal{H}_\alpha(W_t,t)$ is  a martingale
for $\alpha>0$ with respect to the natural filtration.
\end{theorem}

It is   known (see formula \eqref{mo}) that  $\frac{e^{-yW_t}}{\mathbf{E}e^{-yW_t}} = e^{-yW_t- \frac{1}{2}y^2 t}$ is a martingale  and
$$
\varphi(y,u)=\frac{e^{-uy}}{\mathbf{E}e^{-uW_t}}f_{W_t}(y)={\frac {1}{\sqrt {2\pi t}}}e^{-\frac{1}{2t}(y-tu)^2}\in\mathcal{S}(\mathbb{R}\times\mathbb{R}_+).
$$
By Theorem \ref{t02} we have  $\mathcal{H}_{\alpha}(x,t)=\mathcal{A}_{W_t} \{y^{\alpha}\}(x)$.

Generalised function $g{=}\mathcal {L}^{-1}\{y^{\alpha}\}{=}({}^C\mathcal{D}_{0+}^{\alpha}\delta)(x){\in} \mathcal{S}'$ for $\alpha>0$ and function
$e^{-\frac{u^2}{2}t}({}^C\mathcal{D}_{0+}^{\alpha}\delta)(u){\in} \mathcal{D}'^+_a$.
Therefore, by Theorem \ref{t01} $\mathcal{H}_{\alpha}(W_t,t)=\mathcal{A}_{W_t} \{y^{\alpha}\}(W_t)$ is also a martingale.

\begin{theorem} Let $\alpha>0$ and $s<t$, then
\begin{eqnarray*}
{\rm Cov}( \mathcal{H}_\alpha (W_t,t), \mathcal{H}_\alpha (W_s,s))&=&s^\alpha\\
{\rm Corr} ( \mathcal{H}_\alpha (W_t,t), \mathcal{H}_\alpha (W_s,s))&=&\left(\frac{s}{t}\right)^\frac{\alpha}{2}.
\end{eqnarray*}	
\end{theorem}
\begin{proof}
We know that for $t>s$
$${\rm Var}(W_t) = \mathbf{E} W_t^2 = t,$$
$${\rm corr}(W_t, W_s) = \sqrt{\frac{s}{t}} =\rho,\qquad {\rm Cov}(W_t, W_s) = s.$$
Since $\mathcal{H}_\alpha (W_t,t)$ is a martingale $ \mathbf{E}\mathcal{H}_\alpha (W_t,t)= \mathcal{H}_\alpha (W_0,0) =0$
for $\alpha>0$.
\\
Now let us calculate the covariance.
\begin{eqnarray*}
& &{\rm Cov}( \mathcal{H}_\alpha (W_t,t), \mathcal{H}_\alpha (W_s,s)) = \mathbf{E} [\mathcal{H}_\alpha (W_t,t) \mathcal{H}_\alpha (W_s,s)]\\
&=&\frac{1}{2 \pi \sqrt{ts(1{-}\rho^2)} }\int\limits_{-\infty}^\infty \int\limits_{-\infty}^\infty  \mathcal{H}_\alpha (x,t)\mathcal{H}_\alpha (y,s) e^{ -\frac{1}{2(1-\rho^2)}\left(\frac{x^2}{t}-\frac{2\rho xy}{\sqrt{ts}}+\frac{y^2}{s}\right)}dx dy\\
&=&	\frac{1}{2 \pi \sqrt{ts(1{-}\rho^2)} }\int\limits_{-\infty}^\infty \int\limits_{-\infty}^\infty  \mathcal{H}_\alpha (-x,t)\mathcal{H}_\alpha (-y,s) e^{ -\frac{1}{2(1-\rho^2)}\left(\frac{x^2}{t}-\frac{2\rho xy}{\sqrt{ts}}+\frac{y^2}{s}\right)}dx dy.
\end{eqnarray*}
By  \eqref{MT} and formula (5.9) from \cite{SamkKilbMari1993} we obtain
$$
\mathcal{H}_\alpha(-x,t)=\int\limits_{0}^\infty (\,^C\mathcal{D}_{-}^{\alpha}\delta)(u) e^{-ux-\frac {u^2}{2}t} du,
$$
$$
\mathcal{H}_\alpha(-y,s)=\int\limits_{0}^\infty (\,^C\mathcal{D}_{-}^{\alpha}\delta)(v) e^{-vy-\frac {v^2}{2}s} dv,
$$
where $\,^C\mathcal{D}_{-}^{\alpha}$ is given by \eqref{RLDC02},
and
\begin{eqnarray*}
& &
{\rm Cov}( \mathcal{H}_\alpha (W_t,t), \mathcal{H}_\alpha (W_s,s)) =\\
&=&\frac{1}{2 \pi \sqrt{ts(1- \rho^2)} }\int\limits_{0}^\infty\int\limits_{0}^\infty(\,^C\mathcal{D}_{-}^{\alpha}\delta)(u)
(\,^C\mathcal{D}_{-}^{\alpha}\delta)(v)
e^{ -  \frac {u^2}{2}t-\frac {v^2}{2}s}dudv\times\\
&\times&\int\limits_{-\infty}^\infty \int\limits_{-\infty}^\infty e^{ -\frac{1}{2(1-\rho^2)}\left(\frac{x^2}{t}-\frac{2\rho xy}{\sqrt{ts}}+\frac{y^2}{s}\right)-ux-vy}dx dy.
\end{eqnarray*}
For inner integral we have
\begin{eqnarray*}
I(u,v)
&=&\frac{1}{2 \pi \sqrt{ts(1- \rho^2)} } \int\limits_{-\infty}^\infty \int\limits_{-\infty}^\infty e^{ -\frac{1}{2(1-\rho^2)}\left(\frac{x^2}{t}-\frac{2\rho xy}{\sqrt{ts}}+\frac{y^2}{s}\right)-ux-vy}dx dy \\
&=&e^{\frac{1}{2}\left(u^2 t + 2 uv \rho \sqrt{ts} + v^2 s\right)}.
\end{eqnarray*}
$I(u,v)$ is  a moment generating function of a bivariate normal distribution.
\\
Therefore, using formulas (5.16') and (5.20) from \cite{SamkKilbMari1993},  we obtain
\begin{eqnarray*}
& &{\rm Cov}( \mathcal{H}_\alpha (W_t,t), \mathcal{H}_\alpha (W_s,s)) =\\
&=&\int\limits_{0}^\infty\int\limits_{0}^\infty(\,^C\mathcal{D}_{-}^{\alpha}\delta)(u)
(\,^C\mathcal{D}_{-}^{\alpha}\delta)(v)
e^{uv\rho \sqrt{ts} }dudv\\
&=& \int\limits_{0}^\infty\int\limits_{0}^\infty(\,^C\mathcal{D}_{-}^{\alpha}\delta)(u)
(\,^C\mathcal{D}_{-}^{\alpha}\delta)(v)
e^{uvs }dudv
= s^\alpha
\end{eqnarray*}
Subsequently,
\begin{equation*}
{\rm Corr}( \mathcal{H}_\alpha (W_t,t), \mathcal{H}_\alpha (W_s,s)) = \frac{{\rm Cov}( \mathcal{H}_\alpha (W_t,t), \mathcal{H}_\alpha (W_s,s))}{\sqrt{{\rm Var}(\mathcal{H}_\alpha (W_t,t))}{\sqrt{{\rm Var}(\mathcal{H}_\alpha (W_s,s))}}} = \left(\frac{s}{t}\right)^{\frac{\alpha}{2}}.
\end{equation*}
\end{proof}

It is worth noting that when $\alpha=1$, i.e. when  $\mathcal{H}_\alpha (W_t,t)$ coincides with $W_t$, we obtain the standard results for the correlation and covariance of the standard Wiener process, as one would naturally anticipate.
\begin{theorem}
Let $c>0$, then $\mathcal{H}_\alpha (W_t,t)$ is a self-similar process of order $\alpha/2$, i.e. $\left(\mathcal{H}_\alpha(W_{ct}, ct), t \geq 0\right) \stackrel{d}{=} \left(c^{\frac{\alpha}{2}}\mathcal{H}_\alpha(W_{t},t), t \geq 0 \right)$.
\end{theorem}

We have
\begin{eqnarray*}
\mathcal{H}_\alpha(W_{ct}, ct)&=&\mathcal{A}_{W_{c t}}\{ y^{\alpha}\}(W_{ct})=\int\limits_{0}^\infty \mathcal{L}^{-1}\{y^\alpha\}(u)e^{-uW_{ct}-\frac {u^2}{2}c t}du\\&\stackrel{d}{\backsim}&
\int\limits_{0}^\infty \mathcal{L}^{-1}\{y^\alpha \}(u)e^{-u\sqrt{c}W_{t}-\frac {(u \sqrt{c})^2}{2} t}du\left|_{\{\sqrt{c}u=v\}}\right.\\
&=& \frac{1}{\sqrt{c}}\int\limits_{0}^\infty \mathcal{L}^{-1}\{y^\alpha\}\left(\frac{v}{\sqrt{c}}\right) e^{-vW_{t}-\frac {v^2}{2}t}dv.
\end{eqnarray*}
Since
$$
{\mathcal {L}}^{-1}\{y^\alpha\}\left(\frac{v}{\sqrt{c}}\right)=	({}^C\mathcal{D}_{0+}^{\alpha}\delta)\left(\frac{v}{\sqrt{c}}\right)=c^{\frac{\alpha+1}{2}}({}^C\mathcal{D}_{0+}^{\alpha}\delta)\left(v\right),
$$
then
$$
\left(\mathcal{H}_\alpha(W_{ct}, ct), t \geq 0\right) \stackrel{d}{=} \left(c^{\frac{\alpha}{2}}\mathcal{H}_\alpha(W_{t},t), t \geq 0 \right).
$$

\begin{remark}Note that $\mathcal{H}_\alpha(W_{t},t)$ is an example of a non-Gaussian self-similar process, which, informally speaking, is a fractional power of a non-fractional process. It would be interesting to compare it to other known non-Gaussian self-similar processes, such as those discussed in \cite{Tudo2023}, where, again informally speaking, a non-fractional power of a fractional process results in a non-Gaussian self-similar process.
	\end{remark}
\section{Extended Hermite Function}
\subsection{Definition of the Extended Hermite Function}
In this subsection, we take two specific power-normalised parabolic cylinder functions from two different orthogonal sets on the half-line and join them smoothly to create a new function. Our goal is to later use this function to construct  an orthogonal set of functions over the entire real line,  similar to how Hermite polynomials form an orthogonal set. 
\begin{definition} Let $y>0$,and $\alpha >0$ be fixed real parameters.  Let $\mu$, $0<\mu<1$ be a solution of the equation (with fixed $\alpha$):
	\begin{equation}
		\label{Def_mu}
		\frac{1}{\sqrt{\mu}}\frac{\Gamma\left(-\frac{\mu \alpha}{2}\right)}{\Gamma\left(\frac{1-\mu \alpha }{2}\right)} = -  \frac{\Gamma\left(-\frac{\alpha}{2}\right)}{\Gamma\left(\frac{1-\alpha}{2}\right) }.
	\end{equation}
	\label{Extended_Hermite_function_def} 
	We define {\bf the extended Hermite function of order $\alpha$}, denoted as $\tilde{\mathcal{H}}_{\alpha}(x,y)$, by:
	\begin{equation}
		\label{Def_Hermite_Function}
		\tilde{\mathcal{H}}_{\alpha}(x,y)= 
		\begin{cases}
			\mathcal{H}_{\alpha}(x,y), & \text{if  } x \geq 0,\\
			\kappa \mathcal{H}_{\beta}(-x, \mu y), & \text{if  } x<0,
		\end{cases}
	\end{equation}
	where 	$\beta = \mu \alpha$, $\mathcal{H}_{\alpha}(x,y)$ and 	$\mathcal{H}_{\beta}(-x, \mu y)$ are  power-normalised parabolic cylinder functions (see definition \ref{Halpha01}), and the parameter $\kappa$ is given by:
	\begin{equation}
		\label{kappa}
		\kappa = \mu^{-\frac{\mu \alpha}{2}}(2y)^{\frac{\alpha(1- \mu)}{2}}   \frac{\Gamma\left( \frac{1 - \mu \alpha }{2}\right)}{\Gamma\left( \frac{1 - \alpha}{2}\right)}.
	\end{equation}
\end{definition}
\begin{lemma} Let $\alpha >0$. The extended Hermite functions are continuous smooth functions in $x$, $\tilde{\mathcal{H}_\alpha}
	\in C^2 (\mathbb{R}).$
\end{lemma}
\proof
Since the power-normalised parabolic cylinder functions in (\ref{Def_Hermite_Function}) are solutions to the boundary problem (\ref{BVP01}), to prove the lemma it is sufficient to verify that the extended Hermite function is smooth at $0$. Using (\ref{xzero})  it is straightforward  to see that 
\begin{equation*}
\kappa = \frac{\mathcal{H}_{\alpha}(0,y)}{\mathcal{H}_{\beta} (0,\mu y)}=
\mu^{-\frac{\beta}{2}}(2y)^{\frac{\alpha-\beta}{2}}   \frac{\Gamma\left( \frac{1 - \beta }{2}\right)}{\Gamma\left( \frac{1 - \alpha}{2}\right)},
\end{equation*}
i.e. $\mathcal{H}_{\alpha}(0,y) =
\kappa \mathcal{H}_{\beta}(0,\mu y)$ . Moreover,
\begin{eqnarray*}
\left. \frac{\partial}{\partial x} \mathcal{H}_{\alpha}(x,y)\right|_{x=0} &=&
-\frac{\sqrt{\pi } 2^{\frac{\alpha +1}{2}} y^{\frac{\alpha-1}{2}}}{\Gamma \left(-\frac{\alpha }{2}\right)}= 
\mathcal{H}_{\alpha}(0,y) \sqrt{\frac{2}{y}}\frac{\Gamma\left(\frac{1-\alpha}{2}\right)}{\Gamma \left(-\frac{\alpha }{2}\right)} \label{1}
\\ 
&\stackrel{(\ref{Def_mu})}{=} &\kappa
\mathcal{H}_{\beta}(0,\mu y) \frac{1}{\sqrt{\mu}}\sqrt{\frac{2}{y}}\frac{\Gamma\left(\frac{1-\mu \alpha}{2}\right)}{\Gamma \left(-\frac{\mu \alpha }{2}\right)} =\left. \frac{\partial}{\partial x} 
\kappa \mathcal{H}_{\beta}(-x,\mu y)  \right|_{x=0}.
\end{eqnarray*} 
and by (\ref{DxH}) and (\ref{xzero}) we have
\begin{eqnarray*}
	\left. \frac{\partial^2}{\partial x^2} \mathcal{H}_{\alpha}(x,y)\right|_{x=0} &=& \alpha (\alpha-1)\mathcal{H}_{(\alpha-2)}(0,y)= \alpha (\alpha-1)
	\frac{\sqrt{\pi } (2y)^{\frac{\alpha }{2} -1 }}{\Gamma \left(1+\frac{1-\alpha }{2}\right)}\\&= &
-\frac{\alpha}{y} \mathcal{H}_\alpha(0,y)= -\frac{\beta}{\mu y} \kappa \mathcal{H}_\beta(0,\mu y)= 
	\\ 
	&\stackrel{}{=} &\kappa \beta (\beta-1)\mathcal{H}_{(\beta-2)}(0,\mu y) =\left. \frac{\partial^2}{\partial x^2} 
	\kappa \mathcal{H}_{\beta}(-x,\mu y)  \right|_{x=0}. 
 \end{eqnarray*}

\subsection{Orthogonal Set of Extended Hermite Functions}
Before we formulate the theorem regarding the orthogonal set of extended Hermite functions, let us specify the two sets of power-normalized parabolic cylinder functions, orthogonal on the half-line, from which we will construct the orthogonal set of extended Hermite functions on the entire line. Specifically,  we need to determine the sets of real positive numbers $\{\alpha_k\}_{k \in \mathbb{N}}$ and  $\{\beta_k\}_{k \in \mathbb{N}}$ such that  $\{\mathcal{H}_{{\alpha}_k}(x, {y_1})\}_{k \in \mathbb{N}}$  and $\{\mathcal{H}_{{\beta}_k}(x, y_2\}_{k \in \mathbb{N}}$, with fixed $y_1>0, y_2>0$, form  two sets of orthogonal functions on the half-line. According to Theorem \ref{t00}, two real numbers $c_1$ and $c_2$ are sufficient, and equation (\ref{Eq01}) (with $c_1$ and $c_2$ instead of $c$) will provide the desired sets of real positive numbers.
\begin{proposition}  Let $\alpha>0$ be a fixed constant. Suppose  $c =\frac{\Gamma\left(-\frac{\alpha}{2}\right)}{\Gamma\left(\frac{1-\alpha}{2}\right)}$, and let $\mu$, where $0<\mu<1$, be the solution of the equation with fixed $\alpha$:
		\begin{equation}
			\label{unique_root}
		\frac{1}{\sqrt{\mu}}\frac{\Gamma\left(-\frac{\mu \alpha}{2}\right)}{\Gamma\left(\frac{1-\mu \alpha }{2}\right)} = -  \frac{\Gamma\left(-\frac{\alpha}{2}\right)}{\Gamma\left(\frac{1-\alpha}{2}\right) }.
	\end{equation}
	Let $\{\alpha_k\}_{k \in \mathbb{N}}$, where $2k-2 < \alpha_k< 2k$, be the sequence of solutions of the equation
	\begin{equation*} 
		\frac{\Gamma \left(-\frac{x}{2}\right)}{\Gamma \left(\frac{1-x}{2}\right)}  = c.
	\end{equation*}
	Suppose $\beta_k= \mu \alpha_k$; then $2k-2 < \beta_k< 2k$ , and $\beta_k$ is a solution to the equation
	\begin{equation} 
		\label{root_mu_c}
		\frac{\Gamma \left(-\frac{x}{2}\right)}{\Gamma \left(\frac{1-x}{2}\right)}  = - \sqrt{\mu}c.
	\end{equation}
\end{proposition}

 We have $2k-2 < \alpha_k <2k$. Therefore, $(2k-2)<\mu(2k-2) < \mu \alpha_k< \mu 2k<2k$ as $0<\mu<1$. 
\\Define  function $f=f(x,\mu)$ with parameter $\mu>0$  as 
$f(x,\mu)
 = \frac{1}{\sqrt{\mu}} 
 \frac{\Gamma \left( -\frac{\mu x}{2} \right)}
 {\Gamma \left( \frac{1-\mu x}{2} \right)} $. 
 Function $f(x,\mu)$ is continuous and increasing function of $x$ on $[\frac{2k-2}{\mu}, \frac{2k}{\mu}]$. Moreover,  on  $\left(2k-2, 2k\right)$ $\lim_{x\rightarrow2k+1 \pm 1}   f(x,1)  = \pm \infty$. Therefore, equation (\ref{unique_root}) rewrtitten as 
 $$f(x,\mu)= -f(x,1)$$
 has a unique solution on each interval $[2k-2, 2k]$, and equation (\ref{root_mu_c}) has a unique solution  on $[2k-2, 2k]$    $\beta_k= \mu \alpha_k$.

\begin{theorem}  	The extended Hermite functions \label{Extended_Hermite_function_Orthogonality}
	$\{\tilde{\mathcal{H}}_{\alpha_k}(x,y)\}_{k\in\mathbb{N}}$ with fixed parameter $y>0$ form an orthogonal complete set with respect to $x$ on  $\mathbb{R}$, with
	weight function $w(x,y) = \frac{1}{\sqrt{2 \pi y}}e^{-\frac{x^2}{2y}}$:
	\begin{eqnarray*} \label{Eq_0002}
		\langle \tilde{\mathcal{H}}_{\alpha_k},\tilde{\mathcal{H}}_{\alpha_m}\rangle_w&=&	 \frac{1}{\sqrt{2 \pi y}}\int\limits_{-\infty}^\infty \tilde{\mathcal{H}}_{\alpha_k}(x,y)\tilde{\mathcal{H}}_{\alpha_m}(x,y)e^{-\frac{x^2}{2y}}dx=0,\qquad \alpha_k\neq\alpha_m,\nonumber 
		\\||\tilde{\mathcal{H}}_{\alpha_k}(x,y)||_w &=&
	\frac{1}{\sqrt{2} \pi y} \int\limits_{-\infty}^\infty\tilde{\mathcal{H}}_{\alpha_k}^2(x,y)e^{-\frac{x^2} {2y}}dx =y^{\alpha_k}\frac{\left(\psi\left(\frac{1-\alpha_k }{2}\right)-\psi\left(-\frac{\alpha_k
			}{2}\right)\right)}{2\Gamma (-\alpha_k )} \\& +&(\mu y)^{\mu \alpha_k}\frac{\left(\psi\left(\frac{1-\mu \alpha_k }{2}\right)-\psi\left(-\frac{\mu \alpha_k
			}{2}\right)\right)}{2\Gamma (-\mu \alpha_k )}.
	\end{eqnarray*}	
	\end{theorem}

  Let $\beta_k = \mu \alpha_k$.
For $\alpha_k\neq\alpha_m$, we have
\begin{eqnarray*} 
	\langle \tilde{\mathcal{H}}_{\alpha_k},\tilde{\mathcal{H}}_{\alpha_m}\rangle_w&=&	\frac{1}{\sqrt{2 \pi y}}\int\limits_{-\infty}^\infty \tilde{\mathcal{H}}_{\alpha_k}(x,y)\tilde{\mathcal{H}}_{\alpha_m}(x,y)e^{-\frac{x^2}{2y}}dx,\\
	&=& \frac{1}{\sqrt{2 \pi y}}\int\limits_{-\infty}^0 \kappa^2 {\mathcal{H}}_{\beta_k}(-x, \mu y){\mathcal{H}}_{\beta_m}(-x,\mu y)e^{-\frac{x^2}{2y}}dx \\ & &+ \frac{1}{\sqrt{2 \pi y}}\int\limits_{0}^\infty {\mathcal{H}}_{\alpha_k}(x,y){\mathcal{H}}_{\alpha_m}(x,y)e^{-\frac{x^2}{2y}}dx=0.
\end{eqnarray*}
\begin{eqnarray*}
||\tilde{\mathcal{H}}_{\alpha_k}(x,y)||_w &=& 
\frac{1}{\sqrt{2 \pi y}}  \int\limits_{0}^\infty \mathcal{H}_{\alpha_k}^2(x,y)e^{-\frac{x^2} {2y}}dx
+ \frac{1}{\sqrt{2\pi y}} \int\limits_{-\infty}^0 \kappa^2  \mathcal{H}_{\beta_k}^2(-x,\mu y) 
e^{-\frac{x^2}{2y} }dx 
	\\	
	&=& y^{ \alpha_k}   \frac{ \psi\left( \frac{1-\alpha_k }{2}\right)-\psi\left(-\frac{\alpha_k}{2}\right) }{2\Gamma (-\alpha_k )}
	  + (\mu y)^{\mu \alpha_k}\frac{\psi\left(\frac{1-\mu \alpha_k }{2}\right)-\psi\left(-\frac{\mu \alpha_k
		}{2}\right)}{2\Gamma (-\mu \alpha_k )}.
\end{eqnarray*}
The proof of completness is idential to half-line case.

It is well known that the orthogonal decomposition on the real line in terms of Hermite polynomials is related to quantum harmonic oscillator. Similarly,  the asymmetric harmonic oscillator considered in \cite{Chad2023}, is related to the orthogonal decomposition we  above.

\section{Fractional Wiener Chaos Expansion via Extended Hermite Functions}

The question of how to define a fractional It\^o integral is not easy.  One of the way to do it, is to use a power normalised parabolic cylinder function. 
\\
Namely, following Definition \eqref{integral1} we extend the $n$-fold It\^{o} integral of function $g$ as following:
\begin{definition}
	\label{integral_fractional_g}
	Let $\alpha>0$. The fractional $\alpha$-fold It\^{o} integral of $g\in L^2([0,T])$ is defined as
	$$
	I^\alpha g: =\mathcal{H}_\alpha(\Theta,||g||),
	$$
	where
	$$
	\Theta=\int\limits_0^T g(t) dW_t,
	$$
	and $\{W_t\}_{t\geq0}$ is a Wiener process.
\end{definition}
Note, that when $g(t) \equiv 1$ we have
$$
I^\alpha 1: =\mathcal{H}_\alpha(W_t,t),
$$
i.e. the power-normalised parabolic cylinder function is a fractional integral of integrand 1. We will not address the question of whether the extended Hermite function can be considered a fractional integral in this paper; instead, we will reserve the discussion for the second part of this work, to be published separately. It is worth noting  that, as folows from Theorems \ref{t01} and \ref{t02}, $\mathcal{H}_\alpha(W_t,t)$ is a martingale, a property that does not seem to hold for the extended Hermite function $\tilde{\mathcal{H}}_\alpha(W_t,t)$. Further exploration of this topic will be deferred for future publications.

As a further step towards fractional chaos we take the fractional analogue of the building block of the  polynomial chaos defined  in  Definition \ref{integral_tensorproduct}, i.e. the tensor product of extended Hermite functions. Please note that unlike non-fractional case, this building block of the discrete chaos is not a fractional integral, but serves its role in the orthogonal decomposition. 

Suppose $g_1$, $g_2, \dots $ are orthogonal functions in $L^2([0,T])$, and let $\{\alpha_k\}$ be the sequence from Theorem \ref{Extended_Hermite_function_Orthogonality} of real positive non-integer numbers such that $\{\tilde{\mathcal{H}}_{\alpha_k}(x,y)\}_{k\in \mathbb{N}}$ forms a set of orthogonal functions. Then, the analogue of the
Fourier-Hermite function from the original paper by Cameron and Martin,  \cite{CameMart1947}, or the $n^{th}$ fractional polynomial chaos is given by
\begin{equation}
	\prod_{k=1}^{r}\tilde{\mathcal{H}}_{\alpha_{j_k}} \!\!\left( \int\limits_0^T g_{i_k}(t) dW_t, \|g_{i_k}\|^2\right),
\end{equation}
where $r$ is a natural number such that $\{j_k\}_{k=1}^r $ is a set of $r$ pairwise distinct natural numbers such that  $j_1 + \dots+ j_r =n$,  $\{i_k\}_{k=1}^r$ is also a set of $r$ pairwise distinct natural numbers.

Now we are ready to formulate the fractional analogue of Wiener-It\^{o} decomposition. 

Consider the renormalisation of $\int_0^T g_k(t)dW_t$  such that the norms of functions $g_k$  are equal to 1, $\|g_k\|=1$ for all $k$. 
Theorem 1.10 in \cite{DiNuOksePros2009} or Theorem 2.2.4 in \cite{HoldOkseUboeZhan2010}  
can be reformulated in our context as follows (see also \cite{GhanSpan1991}):
\begin{theorem}
	\label{WienerItoExpansionFrac}
	Suppose $F$  is  a  square integrable random variable in $L^2(\mathbb{R}, v(dx))$ with $v(dx) = \frac{1}{\sqrt{2 \pi }}e^{-\frac{x^2}{2}}dx$, and let  $\{ \xi_i\}_{i=1}^\infty$  be a set of independent normally distributed random variables with mean zero and variance one, $\xi_i \sim N(0,1)$.  Let $\{\alpha_k\}_{k=1}^\infty$ be the sequence of real positive non-integer numbers such that $\{\tilde{\mathcal{H}}_{\alpha_k}(x,y)\}_{k=1}^\infty$ forms a set of orthogonal extended Hermite functions. 
	Then there exists a unique representation
	\begin{eqnarray}
		F &=&
		{\bf E} (F)+\sum_{n=1}^\infty \,\, \sum_{j_1\!+\dots\!+j_r\!=n} \,\,\sum_{\substack{i_1, i_2, \dots,i_r;\\ 
				i_l \neq i_m\mbox{\small{for} } l\neq m}} c_{i_1 i_2 \dots i_r}^{j_1 j_2 \dots j_r} \prod_{k=1}^r \tilde{\mathcal{H}}_{\alpha_{j_k}}(\xi_{i_k}, 1)
	\end{eqnarray}
	where $j_k $  and $i_k$ are natural numbers, $c_{i_1 i_2 \dots i_r}^{j_1 j_2 \dots j_r}$ are some constants. The convergence is in $L^2(\mathbb{R},v(dx))$ .
\end{theorem}
The proof is absolutely analogues to the non-fractional case. See for example Theorem 2.6 in \cite{Jans1997}.

It should be noted, that using Theorem 2.11 from \cite{Jans1997} one can extend the results of Theorem \ref{WienerItoExpansionFrac} from $L^2$ to $L^p$, $0 \leq p < \infty$.
\section{Conclusion}
We introduced the function $\mathcal{H}_\alpha(x,y)$, which we refer to as a power-normalized parabolic cylinder function. This function appears to serve as a fractional generalization of the Hermite polynomial on the half-line, as confirmed by our examination of several fundamental deterministic and stochastic properties. In fact, in the same way that the Hermite polynomial $H_n(W_t,t)$ plays the role of ``power $n$ of the Wiener process'', i.e., ``$(W_t)^n$'', the power-normalized parabolic cylinder function $\mathcal{H}_\alpha(W_t,t)$ plays the role of ``the fractional power function of the Wiener process'', i.e., ``$(W_t)^\alpha$".

The power-normalized parabolic cylinder function $\mathcal{H}_\alpha(W_t,t)$ is a martingale and a self-similar non-Gaussian process. Moreover, we propose that the power-normalized parabolic cylinder function can be interpreted as a fractional It\^{o} integral with an integrand of 1, drawing parallels with the non-fractional case.

As a foundational element for the fractional analogue of Wiener chaos, we introduce the extended Hermite function $\tilde{\mathcal{H}}_\alpha(x,y)$ by smoothly joining two power-normalized parabolic cylinder functions, and construct an orthogonal basis on the entire line, $\{\tilde{\mathcal{H}}_\alpha(x,y)\}_{k \in \mathbb{N}}$. Subsequently, from the tensor product of extended Hermite functions as a building block, we construct the fractional analogue of polynomial Wiener chaos.

This paper is the first in a sequence of two papers on Wiener chaos expansion. In the subsequent paper, we plan to explore the properties of the extended Hermite function and present concrete examples demonstrating how fractional Wiener chaos can be used to solve partial differential equations with stochastic input.
\section{Appendix}
\label{sec:8}
\setcounter{section}{8} \setcounter{equation}{0}
Here we present fractional derivatives and integrals from \cite{SamkKilbMari1993} used in the article and some formulas.
\\
Riemann-Lieuwille  fractional integrals $I^\alpha_{0+}$ and $I^\alpha_{-}$ for $x>0$ and for $\alpha>0$ are, respectively
\begin{equation*}
({I}_{0+}^{\alpha}\varphi)(x)=\frac{1}{\Gamma(\alpha )} \int\limits_0^{x} (x-y)^{\alpha-1}  \varphi(y)dy
\end{equation*}		
and
\begin{equation*}
(I_{-}^{\alpha}\varphi)(x)=\frac{1}{\Gamma (\alpha )}  \int\limits_{x}^{\infty}(y-x)^{\alpha-1}\varphi(y)dy.
\end{equation*}
Let $m=[\alpha ]+1$ for non-integer $\alpha>0$. The Riemann-Liouville fractional derivative   $\mathcal{D}_{-}^{\alpha}$ is
\begin{equation} \label{RLD6}
(\mathcal{D}_{-}^{\alpha}\varphi)(x)=
\left\{ \begin{array}{ll}
(-1)^m	\frac{d^{m}}{dx^{m}}(I^{m-\alpha}_{-}\varphi)(x), & \mbox{if $\alpha\notin \mathbb{N}\cup\{0\}$};\\
(-1)^m \frac{d^m}{dx^m}\varphi(x), & \mbox{if $\alpha=m\in \mathbb{N}$}.\end{array} \right.
\end{equation}
Caputo fractional derivatives ${}^{C}\mathcal{D}^\alpha_{0+}$ and ${}^{C}\mathcal{D}^\alpha_{-}$ for $x>0$  are
\begin{equation} \label{RLDC}
({}^{C}\mathcal{D}^\alpha_{0+}\varphi)(x)=
\left\{ \begin{array}{ll}
(I^{m-\alpha}_{0+}\varphi^{{(m)}})(x), & \mbox{if $\alpha\notin \mathbb{N}\cup\{0\}$};\\
\frac{d^m}{dx^m}\varphi(x), & \mbox{if $\alpha=m\in \mathbb{N}$},\end{array} \right.
\end{equation}
\begin{equation} \label{RLDC02}
({}^{C}\mathcal{D}_{-}^{\alpha}\varphi)(x)=
\left\{ \begin{array}{ll}
(-1)^m	 (I^{m-\alpha}_{-}\varphi^{{(m)}})(x), & \mbox{if $\alpha\notin \mathbb{N}\cup\{0\}$};\\
(-1)^m \frac{d^m}{dx^m}\varphi(x), & \mbox{if $\alpha=m\in \mathbb{N}$}.\end{array} \right.
\end{equation}
In formula \eqref{RLDC} in the case of $\alpha=m\in \mathbb{N}$ and $x=0$ we consider $\frac{d^m}{dx^m}\varphi(x)|_{x=0}$ as a limit when $x\rightarrow+0$.
\begin{proposition} For $\beta>0$ Riemann-Liouville  fractional integrals  $I^\beta_{-}$ of $e^{-\frac{y^2}{2t}}$ and $e^{-xy-\frac{y^2}{2t}}$ by $y$ are, respectively
\begin{equation} \label{LI1}
((I_{-}^{\beta})_ye^{-\frac{y^2}{2t}})(\xi)=
t^{\frac{\beta}{2}}e^{-\frac{\xi^2}{4t}}D_{-\beta}\left(\frac{\xi}{\sqrt{t}}\right),
\end{equation}
\begin{equation} \label{LI2}
((I_{-}^{\beta})_ye^{-xy-\frac{y^2}{2}t})(\xi)=
t^{-\frac{\beta}{2}}	e^{\frac{x^2}{2t}-\frac{(x+t\xi)^2}{4t}}D_{-\beta}\left(\frac{x+t\xi}{\sqrt{t}}\right).
\end{equation}
\end{proposition}	
\proof
Let consider first
\begin{eqnarray*}
((I_{-}^{\beta})_ye^{-\frac{y^2}{2t}})(\xi)&=&\frac{1}{\Gamma (\beta )}  \int\limits_{\xi}^{\infty}(y-\xi)^{\beta-1}e^{-\frac{y^2}{2t}}dy=\left|_{\{y=\sqrt{t}\eta\}}\right.\\
&=&\frac{t^{\frac{\beta}{2}}}{\Gamma (\beta )}  \int\limits_{\frac{\xi}{\sqrt{t}}}^{\infty}\left( \eta-\frac{\xi}{\sqrt{t}}\right)^{\beta-1}e^{-\frac{\eta^2}{2}}d\eta
=\left|_{\left\{\eta-\frac{\xi}{\sqrt{t}}=\tau\right\}}\right.\\
&=&\frac{t^{\frac{\beta}{2}}}{\Gamma (\beta )}
\int\limits_{0}^{\infty}\tau^{\beta-1}e^{-\frac{(\tau+\xi/\sqrt{t})^2}{2}}d\tau.
\end{eqnarray*}
Using the representation \eqref{U} of parabolic cylinder function $D_{-\gamma - \frac{1}{2}}(z)$
we can write for $z=x/\sqrt{t}$ and $\gamma=\beta-\frac{1}{2}$
$$
\int\limits_0^{\infty}\tau^{\beta-1}e^{-\frac{(\tau+\xi/\sqrt{t})^2}{2}}d\tau= \Gamma(\beta)e^{-\frac{\xi^2}{4t}}D_{-\beta}\left(\frac{\xi}{\sqrt{t}}\right),
$$		
that gives \eqref{LI1}.
\\
\\
Now, let us consider the second integral:
\begin{eqnarray*}
((I_{-}^{\beta})_ye^{-xy-\frac{y^2}{2}t})(\xi)&=& 	\frac{1}{\Gamma (\beta )}\int\limits_{\xi}^{\infty}(y-\xi)^{\beta-1}e^{-xy-\frac{y^2}{2}t}dy\\
&=&\frac{1}{\Gamma (\beta )}e^{\frac{x^2}{2t}}\cdot e^{-\frac{x^2}{2t}}\int\limits_{\xi}^{\infty}(y-\xi)^{\beta-1}e^{-xy-\frac{y^2}{2}t}dy\\
&=&\frac{1}{\Gamma (\beta )}e^{\frac{x^2}{2t}}\int\limits_{\xi}^{\infty}e^{-\frac{1}{2t}(x+yt)^2}(y-\xi)^{\beta-1}dy\left|_{\{y-\xi=z\}}\right.\\
&=&\frac{1}{\Gamma (\beta )}e^{\frac{x^2}{2t}}\int\limits_{0}^{\infty}e^{-\frac{1}{2t}(x+t(\xi+z))^2}z^{\beta-1}dz\\
&=&	\frac{1}{\Gamma (\beta )}e^{\frac{x^2}{2t}}\int\limits_{0}^{\infty}e^{-\frac{1}{2}\left(\frac{x+t\xi}{\sqrt{t}}+\sqrt{t}z\right)^2}z^{\beta-1}dz\left|_{\{\sqrt{t}z=v\}} \right.\\
&=&\frac{1}{t^{\frac{\beta}{2}}\Gamma (\beta )}e^{\frac{x^2}{2t}}\int\limits_{0}^{\infty}e^{-\frac{1}{2}\left(\frac{x+t\xi}{\sqrt{t}}+v\right)^2}v^{\beta-1}dv.
\end{eqnarray*}
Using formula \eqref{U} we obtain
$$
\int\limits_{0}^{\infty}e^{-\frac{1}{2}\left(\frac{x+t\xi}{\sqrt{t}}+v\right)^2}v^{\beta-1}dv=
\Gamma (\beta )
e^{-\frac{(x+t\xi)^2}{4t}}D_{-\beta}\left(\frac{x+t\xi}{\sqrt{t}}\right),
$$
that gives \eqref{LI2}.

We would like to extend our deepest gratitude to Professor Vassili Kolokoltsov for inviting us to the Isaac Newton Institute for Mathematical Sciences in Cambridge. His insightful comments and suggestions have been invaluable to the advancement of our research. We are indebted to him for his contributions.\\\\
The authors would like to thank the Isaac Newton Institute for Mathematical Sciences, Cambridge, for support and hospitality during the programme
"Fractional differential equations"\, where work on this paper was undertaken. This work was supported by EPSRC grant no EP/R014604/1.

\end{document}